\documentclass[11pt]{article}
\usepackage{amsmath,amssymb,epsf}
\usepackage[french]{babel}
\usepackage[applemac]{inputenc}
\advance \textwidth 3cm
\advance \textheight 3cm
\usepackage{amsfonts}
\usepackage{latexsym}
\newtheorem{theorem}{Th\'eor\`eme}
\newtheorem{lemme}{Lemme}
\newtheorem{corollaire}{Corollaire}
\newtheorem{definition}{D\'efinition}

\newtheorem{remarque}{Remarque}
{\unskip\nobreak\hfil\penalty50\hskip2em\null\nobreak\hfil%
$\Box$\parfillskip0pt\par\medskip}

\newcommand{\Tr}{\mathrm {Tr}}

\newcommand{\trace}{\mathrm{Tr}}

\title{Trace et valeurs propres extr\^{e}mes d'un produit de matrices de Toeplitz. Le cas singulier.}

\author{ Philippe Rambour\thanks{Universit\'{e} de Paris Sud,
      B\^atiment 425; F-91405
Orsay Cedex;
tel : 01 69 15 57 28 ; fax 01 69 15 60 19
      \mbox{e-mail : philippe.rambour@math.u-psud.fr}}
       \and Abdellatif Seghier\thanks{Universit\'{e} de Paris Sud,
        B\^atiment 425; F-91405
Orsay Cedex;
tel : 01 69 15 60 09 ; fax 01 69 15 72 34
       \mbox{ e-mail : abdelatif.seghier@math.u-psud.fr}}}
\date{}
\begin{document}
\maketitle
  \renewcommand{\abstractname}{Rsum}
     \begin{abstract}
\textbf{ Trace et valeurs propres extr\^{e}mes d'un produit de matrices de Toeplitz. Le cas singulier.}\\
Dans un premier th\'eor\`eme nous donnons un d\'eveloppement asymptotique de la trace de la
matrice $T_{N} (f_{1}) T_{N}^{-1} (f_{2})$ avec 
$ f_{1} (\theta)=\vert 1-e^{i \theta}\vert ^{2\alpha_{1}}c_{1}(e^{i \theta})$ et
 $f_{2}(\theta) =\vert 1-e^{i \theta}\vert ^{2\alpha_{2}}c_{2}(e^{i \theta})$, $c_{1}$ et $c_{2}$ \'etant deux fonctions r\'eguli\`eres sur le tore avec $-\frac{1}{2} < \alpha_{1}, \alpha_{2} < \frac{1}{2}$. 
Ensuite nous \'etudions le cas particulier  $\alpha_{1}>0$  et $\alpha_{2}<0$. Nous obtenons alors l'asymptotique de la trace des puissances 
enti\`eres de $T_{N} (f_{1}) T_{N}^{-1} (f_{2})$ 
et nous en d\'eduisons  la limite lorsque $N$ tend vers l'infini  des valeurs propres de cette matrice ce qui nous permet de donner un principe de grandes d\'eviations pour une famille de formes quadratiques de processus al\'eaoires gaussiens stationnaires. 
           \end{abstract}
          \renewcommand{\abstractname}{Abstract}
          \begin{abstract}
\textbf{ Trace and extreme eigenvalues of a product of truncated Toeplitz matrices. The singular case.}\\
      In a first theorem we give an asymptotic expansion of $\trace \left( T_{N} (f_{1}) T_{N}^{-1} (f_{2})\right)$ where 
      $ f_{1} (\theta)=\vert 1-e^{i \theta}\vert ^{2\alpha_{1}}c_{1}(e^{i \theta})$ and $f_{2}(\theta) =
      \vert 1-e^{i \theta}\vert ^{2\alpha_{2}}c_{2}(e^{i \theta})$, with $c_{1}$ and $c_{2}$ are two regular 
      functions of the torus  and  $-\frac{1}{2} < \alpha_{1}, \alpha_{2} < \frac{1}{2}$.
      In a second part of this work we study the particular case where $\alpha_{1}>0$ and $\alpha_{2}<0$. Then we obtain the asymptotic
of $\trace \left( T_{N} (f_{1}) T_{N}^{-1} (f_{2})\right)^s$ for $s\in \mathbb N^*$ that provides us the limits when $N$ goes to the infinity of the extreme eigenvalues of 
this matrix. This last result allows us to give a large deviation principle for a family of quadratic forms of stationnary 
process.
   \end{abstract}
   
%


\section{Introduction}
Si $f\in L^1 (\mathbb T)$ on appelle matrice de Toeplitz d'ordre $N$ et de symbole $f$, et on note $T_{N}(f)$, la matrice 
d\'efinie par $\left(T_{N}(f)\right)_{k+1,l+1}= \widehat f (l-k)$ pour $0 \le k \le N$ et $0 \le l \le N$, o\`u 
$ \widehat h (s)$ d\'esigne le coefficient de Fourier d'ordre $s$ de la fonction $h$. On dira que le symbole $f$ est r\'egulier si 
la fonction $f$ est strictement positive sur le tore, et que le symbole $f$ est singulier si la fonction $f$ admet des z\'eros ou 
des p\^{o}les sur le tore. Une bonne approche des matrices de Toeplitz peut se trouver dans \cite{Bo.3}.

Un probl\`eme de l'\'etude des matrices de Toeplitz est d'\'etablir la trace du produit de deux matrices de Toeplitz, ou
m\^{e}me d'une puissance d'un produit de matrice de Toepliz.  
Cette \'etude intervient autant dans le domaine de l'analyse (recherche des valeurs propres) que des la probabilit\'es 
( th\'eor\`emes de limite centrale, principes de grandes d\'eviations). Le probl\`eme de l'\'etude de $\trace \left ( T_{N}(f) T_{N}(g) \right)^s $ 
est un grand classique de la litt\'erature consacr\'ee aux matrices de Toeplitz et \`a l'\'etude  des processus al\'eatoires Gaussiens. Il a 
\'et\'e en particulier \'etudi\'e par Avram (\cite{Avram1})dans le cas r\'egulier ou pour une puissance deux, et par Fox et Taqqu (\cite{F.T.1}, \cite{F.T.2}),  dans un cadre plus g\'en\'eral. Il faut aussi citer Grenander et Szeg\^{o} \cite{GS}, Ibragimov \cite{Ib2}, 
Rosenblatt \cite{Rosenb}, Taniguchi \cite{TANI2},
Dalhaus \cite{Dahl}, Giraitis et Surgalis \cite{Gira}, Ginovyan \cite{Ginov0},  Taniguchi and Kakizawa \cite {TANI1}, Ginovyan et Sahakyan \cite{Ginov2}, \cite{Ginov3}
et \cite{Ginov4}, et Lieberman et Philips \cite{LIEBER1}. Un bon r\'esum\'e des principaux acquis peut se trouver dans 
\cite{Ginov1}.

La question qui peut \'egalement se poser est de conna\^{i}tre la trace de $\left(T_{N}(f) T_{N}^{-1} (g)\right)^s$. Une fa\c con de faire 
peut consister \`a se ramener \`a la trace de $\left(T_{N}(f) T_{N} (g^{-1})\right)^s$, mais cette m\'ethode n'est pas toujours 
satisfaisante, surtout dans le cas singulier. Dans cet article nous proposons un d\'eveloppement asymptotique d'ordre 1 ou 2, 
suivant les cas, de $\trace \left( T_{N} (f_{1}) T_{N}^{-1} (f_{2})\right) $ dans le cas o\`u $f_{1} (\theta) = \vert 1 - e^{i \theta} \vert ^{2 \alpha_{1}} c_{1} (\theta) $
et $ f_{2}  (\theta) = \vert 1 - e^{i \theta} \vert ^{2 \alpha_2} c_{2} (\theta) $ avec $- \frac{1}{2} < \alpha_{1}, \alpha_{2} < \frac{1}{2}$
et o\`u $c_{1}$ et $c_{2}$ sont deux fonctions r\'eguli\`eres sur le tore. Les m\'ethodes que nous utilisons sont diff\'erentes de celles de      
 Fox et Taqqu, mais utilisent des travaux ant\'erieurs (voir \cite{RS09} et \cite{RS10}). Nous mettons en \'evidence les diff\'erents cas qui peuvent se pr\'esenter, et la suite de notre travail est consacr\'ee \`a l'\'etude du cas $\alpha_{1} \in ]0, \frac{1}{2}[$ et 
$ \alpha_{2}\in ]- \frac{1}{2} ,0[$ (cas 1) ii)). Nous  donnons dans ce cas une expression asymptotique 
de $\trace \left( T_{N} (f_{1}) T_{N}^{-1} (f_{2})\right)^s $ et aussi de $\trace \left( T_{N} (f_{1}) T_{N} (f_{2}^{-1})\right)^s $ 
(ce qui revient \`a donner  $\trace \left( T_{N} (f_{1}) T_{N}^{-1} (f_{2})\right)^s $ pour $0 <\alpha_{1},\alpha_{2}< \frac{1}{2}$ 
compl\'etant ainsi des r\'esultats de  Taniguchi(\cite {TANI2}, \cite {TANI1} et de Lieberman et Phillips \cite{LIEBER1}. 
Nous \'etudions ensuite  les cons\'equences de ces r\'esultats pour les valeurs propres de $T_{N} (f_{1}) T_{N} (f_{2}^{-1})$ 
et nous en d\'eduisons  un r\'esultat probabiliste. Dans le cas d'une matrice de Toeplitz le comportement des 
valeurs propres ob\'eit \`a certains principes bien connus. On sait par exemple que si $ \lambda_{i}^{(N)}$ sont les valeurs propres d'une matrices de Toeplitz d'un symbole $f$ avec $\inf _{\theta \in \mathbb T} f(\theta) =m$ et $\sup_{\theta \in \mathbb T} f(\theta) =M$
alors $ \displaystyle{\lim_{ N \rightarrow + \infty}\left ( \inf_{1\le i \le N} \lambda_{i}^{(N)}\right)  = m}$ et 
 $ \displaystyle{ \lim_{N \rightarrow + \infty}\left ( \sup_{1\le i \le N} \lambda_{i}^{(N)}\right)  = M}$. Cette propri\'et\'e n'est \'evidemment plus vraie pour un produit de matrices de Toeplitz ni \`a plus forte raison pour le produit d'une matrice de Toeplitz avec une matrice hermitienne. 
 Si on \'etudie une forme quadratique $W_{N}$ d\'efinie par $W_{N} = \frac{1}{N} X^{(N)*} M_{N} X^{(N)}$ avec $(X_{n})$ un processus stationnaire centr\'e gaussien et la matrice $M_{N} $ une matrice hermitienne (voir, entre autre, \cite{Sa.Ka.Ta},  \cite {B.G.R.1} et \cite{B.G.L.1})on ne peut donc pas alors appliquer le th\'eor\`eme de G\^{a}rtner-Ellis (\cite {DemZei}) pour obtenir un principe de grandes d\'eviations pour 
 $W_{N}$. Dans le cas o\`u $M_{N}= T_{N}(f)$ et o\`u $f,g\in L^\infty (\mathbb T)$ ( $g$ \'etant la densit\'e spectrale de 
 $X^{(N)}$) \cite {B.G.R.1} et \cite{B.G.L.1} donnent des solutions \`a ce probl\`eme.
Cela leur permet notamment de donner un principe de grandes d\'eviations  dans l'\'etude du rapport de vraissemblance de deux processus gaussiens  stationnaires (\cite{CO-DA}, \cite{CO-DA1},\cite{DA},\cite{Bouaziz} \cite{Barone}) dans le cas o\`u les densit\'es spectrales sont r\'eguli\`eres. 
 Supposons maintenant  que $X^{(N)}$ admette pour densit\'e spectrale une fonction $f_{1}$
 et que $M_{N} =T_{N}^{-1} (f_{2})$. Dans le cas $\alpha_{1}$positif et $\alpha_{2}$ n\'egatif en notant $\mu_{ i }^{(N)}$ les valeurs propres de 
 $T_{N}(f_{1}) T_{N}^{-1} (f_{2})$ nous obtenons, en \'etudiant la convergence  de la mesure 
 $\displaystyle{ \sum_{i=0}^N \delta _{\mu_{i}^{(N)}} }$, que   
 $ \displaystyle{ \lim_{ N \rightarrow + \infty}\left ( \inf_{1\le i \le N} \mu_{i}^{(N)}\right) } = 
 0$ et $ \displaystyle{ \lim_{N \rightarrow + \infty}\left ( \sup_{1\le i \le N} \mu_{i}^{(N)}\right) } = \Vert \frac{f_{1}}{f_{2}}\Vert _{\infty}$. 
 Nous pouvons alors obtenir sur un intervalle maximal la limite de la suite de fonctions (voir \cite{GS})
\begin{align*}
L_{N}(t)  &=\frac{1}{N} \ln E (e^{N t W_{N}})\\
 &= - \frac{1}{2N} \ln \det \left ( I_{N}- 2 t T_{N}^{1/2} (f_{1})  T_{N}^{-1} (f_{2})T_{N}^{1/2}(f_{1}) \right) 
 = - \frac{1}{2N} \sum_{i=1}^N \ln (1-2 t \mu_{i}^{(N)} ).
 \end{align*}
  Toujours dans le cas 1) ii) avec $\alpha_{1}>0$ et $\alpha_{2} <0$ nous retrouvons un th\'eor\`eme de limite centrale similaire \`a celui donn\'e par Fox et Taqqu dans 
  \cite{F.T.1} mais dans lequel nous pouvons consid\'erer la forme quadratique $X^{(N)*} T_{N} (f_{2})^{-1}X^{(N)}$ au lieu 
  de $X^{(N)*} T_{N} (f_{2}).X^{(N)}$.
  
  Dans un prochain travail nous nous attacherons \`a d\'evelopper les  cas diff\'erents du cas 1) ii) qui interviennent dans le th\'eor\`eme
    \ref{thprincipal}.  Nous aurons deux objectifs : d'abord \'etablir des expressions asymptotiques pour les traces de $\left(  T_{N} (f_{1})
  T_{N} ^{-1} (f_{2})\right)^{s} $ puis utiliser ces expressions pour donner des th\'eor\`emes de limites centrales pour des formes quadratiques 
  du type $X^{(N)*} T_{N} (f_{2})^{-1}X^{(N)}$ ou $X_{N}$ est un processus al\'eatoire centr\'e gaussien stationnaires de densit\'e spectrale $f_{1}$. 
\section{Principaux r\'esultats}
\subsection { Trace de produits de matrices de Toeplitz.}
Dans la suite nous noterons $\chi$ la fonction d\'efinie par $\chi(\theta) =e^{i \theta}$. Pour tous les r\'eels $\nu$ positifs nous consid\'ererons aussi les ensembles 
$A(\mathbb T, \nu) = \{h \in L^2 (\mathbb T) / \displaystyle{ \sum_{k\in \mathbb Z} \vert k \vert ^\nu \hat h(k)} <\infty\}$, 
o\`u $\hat h (k)$ d\'esigne le coefficient d'ordre $k$ de la fonction $h$.
\begin{theorem} \label{thprincipal} 
On consid\`ere deux fonctions $f_{1}$ et $f_{2}$ 
d\'efinies sur le tore $\mathbb T$ par 
$f_{1} = \vert 1-\chi \vert ^{2\alpha_{1} }c_{1}$ et $f_{2} = \vert 1-\chi \vert ^{2\alpha_{2} }c_{2}$ 
o\`u $c_{1}$ et $c_{2}$ sont deux fonctions r\'eguli\`eres  dans $A(\mathbb T, \frac{3}{2})$ et $- \frac{1}{2} < \alpha_{1},\alpha_{2}<\frac{1}{2}$.
On a alors les r\'esultats suivants :
\begin{enumerate}
\item
Si $\alpha_{2} \in ]- \frac{1}{2},0[$
\begin{itemize}
\item [i)]
dans le cas ou  $\frac{1}{2} > \alpha_{2} -\alpha_{1}>0$
  $$\Tr  \left( T_{N} (f_{1}) T_{N}^{-1} (f_{2}) \right) = 
   \trace \left(T_{N} \frac {f_{1}}{f_{2}}\right) + O(N^{2\alpha_{2}-2\alpha_{1}} ),$$
\item [ii)]
dans le cas ou  $ 0> \alpha_{2} -\alpha_{1}> -1$
\begin{eqnarray*}
\Tr  \left( T_{N} (f_{1}) T_{N}^{-1} (f_{2}) \right) &=&
   \trace \left(T_{N} \frac {f_{1}}{f_{2}}\right)  +  \widehat f_{1} (0) \left( \langle \ln f_{2}^{-1}, \frac{1}{f_{2}} \rangle_{2,1/2} \right)
  \\
 &+& C_{1}(f_{1},f_{2}) +
 O(N^{\max (2\alpha_{2}-2\alpha_{1}, -1))}).
 \end{eqnarray*}
 \end{itemize}
 Si les fonctions $f_{1}$ et $f_{2}$ sont des fonctions paires on a 
 $$ C_{1}(f_{1},f_{2}) =-2 (\alpha_{2}+1) \left( \langle  \ln f_{2}^{-1}, \frac{f_{1}}{f_{2}} \rangle_{2,1/2} - \langle  \ln f_{2}^{-1}, \frac{1}{f_{2}} \rangle_{2,1/2} \widehat {f_{1}}(0)
 - \langle  \ln f_{1}, \frac{1}{f_{2}} \rangle_{2,1/2}\right).$$
 \item
 Si $\alpha_2 \in ]0, \frac{1}{2}[$
 \begin{itemize}
 \item [i)] 
 Dans le cas o $ \frac{1}{2} < \alpha_2- \alpha_1 <1 $
 on a 
 $$ \Tr  \left( T_{N} (f_{1}) T_{N}^{-1} (f_{2}) \right) =  N^{2\alpha_2 -2 \alpha_1}  C_{2} (f_{1},f_{2})+
  o(N^{2\alpha_2 -2 \alpha_1})$$
  
    \item [ii)] 
   Dans le cas o $- \frac{1}{2} < \alpha_2- \alpha_1 <\frac{1}{2} $
 on a 
 \begin{itemize} 
 \item [a)]
 Si $\alpha_{1} <0$ 
$$
\trace  \left( T_{N} (f_{1}) T_{N}^{-1} (f_{2}) \right) = \trace \left(T_{N} \frac {f_{1}}{f_{2}}\right)+
N^{2\alpha_2 -2 \alpha_1}  C_{3} (f_{1},f_{2}) + o(N^{2\alpha_2 -2 \alpha_1})
$$

\item [b)]
Si $\alpha_{1} >0$
$$ 
\Tr  \left( T_{N} (f_{1}) T_{N}^{-1} (f_{2}) \right) = \trace \left(T_{N} \frac {f_{1}}{f_{2}}\right)
+ 
\left( \int_{1/2}^1 \tilde G_{\alpha_{2}} (x) dx \right)
  \frac{2  N^{2\alpha_{2} } \widehat f_{1} (0)}{c_{2} (1) \Gamma ^2 (\alpha_{2})}  + 
  o(N^{2\alpha_{2}})
  $$
 avec 
  $$ \tilde G_{\alpha_{2}} (x) = \int_{0}^x  t^{2\alpha_{2}-2}\left( (1-t)^{2\alpha_{2}}-1 \right) - 
  t^{2\alpha_{2}}(1-t)^{2\alpha_{2}-2} dt dx + \frac{x^{2\alpha_{2}-1}} {(2\alpha_{2}-1)},$$
  et $x \in [1/2,1]$.

  \end{itemize}
\end{itemize}
\end{enumerate}
\end{theorem}
\begin{remarque}
Les constantes $C_{i}(f_{1},f_{2})$ $1\le i \le 3$, peuvent \^{e}tre obtenues \`a partir de la d\'emonstration du 
th\'eor\`eme.
\end{remarque}

Dans l'\'enonc\'e suivant nous \'etudionsun cas particulier du th\'eor\`eme \ref{thprincipal}.
et nous noterons $f_{2,t}$ la fonction 
$$ f_{2,t}(\theta) =  t f_{1}(\theta) + f_{2} (\theta),$$
 \begin{theorem}\label{thpuissance}
 Si $f_{1}$ et $f_{2}$ sont deux fonctions paires d\'efinies sur le tore et v\'erifiant les hypoth\`eses  du th\'eor\`eme \ref{thprincipal} avec de plus $\alpha_{1}$
 positif et $\alpha_{2}$ n\'egatif
 on a, pour tout entier naturel $s$ non nul
$$
\Tr  \left( T_{N} (f_{1}) T_{N}^{-1} (f_{2}) \right)^s -
   \trace \left(T_{N} \left(\frac {f_{1}}{f_{2}}\right)^s\right)  = (-1)^{s-1}\frac {\Psi_{1}^{(s-1)} (0)}{(s-1)!} + o(1).$$
 avec 
   $$ \Psi_{1} (t) =
   \widehat f_{1} (0) \left( \langle \ln f_{2,t}^{-1}, \frac{1}{f_{2,t}} \rangle_{2,1/2} \right)
  - C_{1}(f_{1},f_{2,t}).
 $$
\end{theorem}
\begin{corollaire}\label{corollaire}
On consid\`ere deux fonctions $h_{1}$ et $h_{2}$ paires 
d\'efinies sur le tore $\mathbb T$ par 
$h_{1} = \vert 1-e^{i \theta} \vert ^{2\alpha_{1} }d_{1}$ et $h_{2} = \vert 1-e^{i \theta} \vert ^{2\alpha_{2} }d_{2}$ 
o\`u $d_{1}$ et $d_{2}$ sont deux fonctions r\'eguli\`eres appartenant \`a  $A(\mathbb T, \frac{3}{2})$. 
On suppose de plus 
que $\frac{1}{2}>\alpha_1, \alpha_{2}>0$. Nous pouvons alors crire, 
pour tout entier $s \ge 1$ 
$$ \trace\left (T_N h_1 T_N h_2 \right) ^s = \frac{N}{2 \pi} 
\int_0^{2 \pi} \left((h_1 h_2)(\theta)\right)^s d\theta +o(N).$$
\end{corollaire}
\subsection{Quelques applications aux grandes d\'eviations et aux valeurs propres}
Dans tout ce paragraphe nous consid\'erons encore  
$f_{1} = \vert 1-\chi\vert ^{2\alpha_{1} }c_{1}$ et $f_{2} = \vert 1-\chi \vert ^{2\alpha_{2} }c_{2}$ 
deux fonctions paires avec $\alpha_{1}$ positif et $\alpha_{2}$ n\'egatif et $c_{1}$, $c_{2}$ deux fonctions r\'eguli\`eres dans 
$A(\mathbb T, \frac{3}{2})$.
\begin{lemme} \label{Szego}
Si $\mu_{i}^{(N)}$ $1\le i \le N$ d\'esignent les valeurs propres de $T_{N}(f_{1}) T_{N}^{-1} (f_{2})$, class\'ees dans l'ordre 
croissant, alors la suite de mesures $\displaystyle{ \sum_{i=1}^N \delta _{\mu_{i}^{(N)}} }$ converge au sens faible (ou en loi) vers 
$P_{\frac{f_{1}}{f_{2}}}$ la mesure image de la mesure de Lebesgue sur le tore par $\frac{f_{1}}{f_{2}}$. 
\end{lemme}
Ce lemme admet comme corollaire imm\'ediat les deux th\'eor\`emes suivants 
\begin{theorem} \label{VP}
Avec les hypoth\`eses et notations pr\'ec\'edentes 
$$ \lim_{N \rightarrow + \infty} \mu_{1}^{(N)} =0, \quad  \lim_{N \rightarrow + \infty} \mu_{N}^{(N)} = \Vert \frac {f_{1}}{f_{2}} \Vert _{\infty}.$$
\end{theorem}
On consid\`ere maintenant pour tout entier $N$ la forme quadratique $W_{N}$ d\'efinie par 
$$ W_{N}= \frac{1}{2 N} ^t X^N \left(  T_{N}^{-1} f_{2} \right) X^N$$
o\`u $ X^N$ un processus de densit\'e spectrale 
$f_{1} = \vert 1 - \chi \vert ^{2\alpha_{1}} c_{1}$ et $f_{2}= \vert 1 - \chi \vert ^{2\alpha_{2}} c_{2}$, $c_{1}$ et $c_{2}$ 
\'etant deux fonctions r\'eguli\`eres sur le tore. On consid\`ere alors la suite de fonctions 
$$L_{N} (t) = \frac{1}{2N} \left(  - \frac{1}{2} 
\sum_{i=1}^N \ln (1- 2\mu_{i}^N t)\right)$$
o\`u $ (\mu_{i}^N)_{i=1\cdots N} $ sont les valeurs propres de 
$A_{N} = T_{N}^{1/2} f_{1}T_{N}^{-1} f_{2}T_{N}^{1/2} f_{1}$ qui sont aussi celles 
de $T_{N}(f_{1}) T_{N}^{-1} f_{2}$. Alors nous pouvons \'ecrire, en posant 
$\Delta = ]-\delta^{-1} ,\delta ^{-1}[$ avec $\delta = \Vert \frac{f_{1}}{f_{2}}\Vert_\infty$
\begin{theorem} \label{deviant1}
Si $\frac{f_{1}}{f_{2}} \in L^{\infty} ( \mathbb T)$ avec $\alpha_{1}>0$ et $\alpha_{2}<0$  
 on a pour tout $t$ tel que $2t \in \Delta$ 
$$ L_{N}(t) =- \frac{1}{4\pi}  \int_{0}^{2\pi} \ln \left( 1- 2 t \frac{f_{1}}{f_{2}} (\theta) \right) 
d\theta + o(1),$$
et 
$$ \lim_{N\rightarrow + \infty } \left ( N L_{N}(t) + \frac{N}{4 \pi}\int_{0}^{2\pi} \ln \left( 1- 2 t \frac{f_{1}}{f_{2}} (\theta) \right) 
d\theta \right) = \frac{\Psi (2t)} {2}$$
avec $ \Psi(t) = \sum_{l=1}^\infty \frac{(-1)^{l+1}}{l} (t)^l  \frac{\Psi_{1}^{(l)} (0)}{l!} $
  \end{theorem}
 \begin{remarque}
 Ce th\'eor\`eme revient \`a dire que si  $\frac{f_{1}}{f_{2}} \in L^{\infty} ( \mathbb T)$ et $\alpha_{1}\alpha_{2}<0$ 
 alors  $(W_{N})$ satisfait une SLDP ( Sharp Large Deviation Principle) pour une fonction  $L^*$ qui est  le dual de 
 Fenchel-Legendre est $L(t) = \frac{1}{4\pi} \int_{0}^{2\pi} \ln \left( 1 - 2 t \frac{f_{1}}{f_{2}} (\theta) \right) d\theta$.
  \end{remarque}
  D'autre part en posant $m_{N}= E(^t X^N \left(  T_{N}^{-1} f_{2} \right) X^N)$ nous pouvons \'enonc\'er le th\'eor\`eme 
\begin{theorem} \label{central1}
La variable al\'eatoire :
$$ \frac{^t X^N \left(  T_{N}^{-1} f_{2} \right) X^N - m_{N}}{\sqrt N} $$
converge en loi vers une variablle al\'eatoire qui suit une loi normale centr\'ee de variance 
$1 \pi \int _{-\pi}^{\pi}\left( \frac{f_{1}(\theta)}{f_{2}(\theta)}\right)^2 d \theta.$
\end{theorem}
\begin{remarque}
La d\'emonstration de ce th\'eor\`eme est bien s\^{u}r parfaitement identique \`a celle du th\'eor\`eme du m\^{e}me 
genre donn\'e dans \cite{F.T.1}
\end{remarque}
Pour d\'emontrer le th\'eor\`eme \ref{thprincipal} nous allons devoir, dans un premier temps, donner une expression asymptotique simple
pour $N$ suffisamment grand, des coefficients du polyn\^{o}me pr\'edicteur de degr\'e $N$ d'une fonction  $f $ admettant une singularit\'e
dordre $\alpha$ comprise entre $-\frac{1}{2}$ et $\frac{1}{2}$. C'est ce que nous faisons, apr\`es un bref rappel, dans la partie suivante.
\section{Asymptotique des coefficients du polyn\^ome pr\'edicteur}
\subsection{ Definition et propri\'et\'es fondamentales du polyn\^{o}me pr\'edicteur}
\begin{definition}
Si $h$ une fonction positive dans $L^1 (\mathbb T)$ on appelle polyn\^{o}me pr\'edicteur de degr\'e $N$ de $h$ le 
polyn\^{o}me $K_{N}$ d\'efinie par $\displaystyle{ K_{N}= \sum_{k=0}^N \frac{ (T_{N}^{-1})_{k+1,1}} 
{\sqrt { (T_{N}^{-1})_{1,1}}} z^k}$.
\end{definition}
Nous avons alors les r\'esultats suivants (voir \cite{Ld})
\begin{theorem} \label{theoremmoinsdeux}
Si $h$ une fonction positive dans $L^1 (\mathbb T)$ et $K_{N}$ son polyn\^{o}me pr\'edicteur de degr\'e $N$ alors  
\begin{itemize}
\item
$$ \forall s, \; -N \le s \le N \; \widehat {\vert K_{N}\vert ^{-2}} (s) = \widehat{h} (s).$$
\item
$K_{N}$ ne s'annule pas sur le tore.
\end{itemize}
\end{theorem}
Ce th\'eor\`eme admet la cons\'equence imm\'ediate suivante 
\begin{theorem} \label{theoremmoinstrois}
Avec les hypoth\`eses du th\'eor\`eme pr\'ec\'edent nous avons 
$$ T_{N}\left( \vert K_{N}\vert ^{-2} \right)= T_{N}(h).$$
\end{theorem}
Le calcul de $T_{N}^{-1}(h)$ s'en trouve alors facilit\'e gr\^{a}ce au lemme suivant qui a \'et\'e \'etabli dans \cite{RS07} et 
qui est une version alg\'ebrique de la formule de Gohberg-Semencul \cite{GoSe}. 
\begin{lemme} \label{Gohberg}
Si $P =\displaystyle { \sum_{u=0}^N \delta _{u} z^u }$ un polyn\^{o}me de degr\'e $N$ sans z\'eros sur le tore 
on a 
$$ T_{N}^{-1} \left( \vert P \vert ^{-2} \right) = (\bar \delta _{0}\delta _{l-k}+ \cdots + \bar \delta _{k} \delta _{l}
- ( \delta _{N-k} \bar \delta _{N-l} + \cdots \delta _{N}\delta _{N+k-l}).$$
\end{lemme}
Dans cette partie nous consid\'erons une fonction  $f$  d\'efinie  par $f = \vert1-\chi\vert^{2\alpha} c$ 
o\`u $c$ est une fonction r\'eguli\`ere sur le tore avec $c_{1}=g_{1}\bar g_{1}$, $g_{1} \in H^{2+}$ et
$ g =(1-\chi)^\alpha g_{1}$ . On pose $\beta_{k}^{(\alpha)}= \widehat {g^{-1}} (k)$ et on note par $\beta_{k,N}$
le coefficient de $\chi^k$ du polyn\^{o}me pr\'edicteur de degr\'e $N$ de $f$. D'autre part nous nous pla\c cons dans le cas 
$-\frac{1}{2} <\alpha<\frac{1}{2}$.
\begin{theorem} \label{theoremezero}
On consid\`ere une fonction $f$ comme ci-dessus
Alors si $n_{1}$ est entier fix\'e, ind\'ependamment de $N$,
$$ \beta_{k,N} = \beta_{k}^{(\alpha)} (1-\frac{k}{N})^\alpha \left(1+o(1)\right)$$
pour tout entier $k \in [0, N -n_{1} ]$, uniform\'ement par rapport \`a $N$.
\end{theorem}
\begin{remarque}
Dans la pratique $n_{1}$ est choisi de mani\`ere \`a ce que pour tout entier $u\ge n_{1}$ on ait 
$\beta_{u}^{(\alpha)} = \frac{1}{g_{1}(1)} \frac{u^{\alpha-1}}{\Gamma(\alpha)}$ avec la pr\'ecision n\'ecessaire.
\end{remarque}
\begin{remarque}
Ce th\'eor\`eme peut se lire \\
$ \forall \epsilon >0 \quad  \exists N_{0} \quad \mathrm{t. q.}\quad  \forall N \ge N_{0}
\quad \forall k, \quad 0\le k \le N-n_{1} \quad  \exists R_{k}, \quad \vert R_{k} \vert \le \varepsilon$
 tel que  
$$ \beta_{k,N} = \beta_{k}^{(\alpha)} (1-\frac{k}{N})^\alpha \left(1+R_{k}\right)$$
\end{remarque}

\subsection {D\'emonstration du th\'eor\`eme \ref{theoremezero}}
Pour dmontrer ce thorme nous allons dcouper l'intervalle 
$[0,N \delta]$ en deux,  savoir $[0, \delta_1]$ et 
$[\delta_1,\delta]$. La dmonstration du thorme sur 
$[\delta_1,\delta]$ est assez rapide. Dans 
\cite{RS10} nous avons dmontr
\begin{theorem} \label{theomoinsun}
Si $-\frac{1}{2} < \alpha<\frac{1}{2}$, $\alpha\not=0$ nous avons 
pour $0<x<1$ 
$$ g_1(1) \beta_{[Nx],N} = 
N^{\alpha-1} \frac{1}{\Gamma (\alpha) } x^{\alpha-1}(1-x)^\alpha 
+ o(N^{\alpha-1})$$
uniformment en $x$ dans $[\delta_1,\delta_2]$, pour 
$0< \delta_1 <\delta_2 <1.$
\end{theorem}
Ce thorme est quivalent  
$$ g_1(1) \beta_{k,N} = 
 \frac{1}{\Gamma (\alpha) } k^{\alpha-1}(1-\frac{k}{N} )^\alpha 
+ o(N^{\alpha-1})$$
pour tout entier naturel 
$k$ dans $ [N\delta_1, N \delta_2]$, uniformment par rapport  $N$.
Cette remarque, jointe au rsultat (voir \cite{Zyg}) 
$ \beta^{(\alpha)} _k = \frac{1}{g_1(1) \Gamma (\alpha)} k^{\alpha-1}
\left( 1+o(1)\right)$, uniformment par rapport  $k$ pour $k$ assez grand, permet d'obtenir le rsultat sur $[N \delta_1, N \delta]$.\\
Pour d\'emontrer le th\'eor\`eme dans l'intervalle $[0, N\delta_1]$, il nous faut revisiter
les rsultats de \cite{RS10}.\\
Nous avons obtenu dans cet article le lemme 
\begin{lemme} \label{lemmemoinsun}
On suppose $\alpha \in ]-\frac{1}{2}, \frac{1}{2}[.$ Alors, sous les hypothses du thorme \ref{theomoinsun} nous pouvons crire pour $N$ suffisamment grand et pour 
$ 0 \le k \le \delta N $, $0 <\delta <1$ 
$$ 
\beta_{k,N}
= \beta_k^{(\alpha)} - \frac{1}{N} \sum_{u=0}^k \beta_{k-u}^{(\alpha)}
\left( F_\alpha \left( \frac{u}{N} \right) \right) \left(1+o(1)\right).
$$
la fonction $z \rightarrow F_\alpha (z)$ est continue et drivable 
sur tout compact de $[0,1[$ et de plus pour tout rel dans $[0, \delta]$ on a 
$$ \vert F_\alpha (z) \vert \le K_0 (1 + \vert \ln (1-z)\vert ),$$
o $K_0$ est une constante indpendante de $N$.
\end{lemme}
\begin{remarque} 
En utilisant le d\'eterminant de  $T_N ( \vert 1-\chi \vert ^{2 \alpha} f_1)$ et 
la formule d'Hartwig-Fisher nous obtenons facilement  
$$ F_{\alpha}(0) = \alpha^2 +o(1).$$
\end{remarque}
Si
 $k\in [0, k_{0}]$ et $N$ suffisamment grand le lemme \ref{lemmemoinsun} permet d`\'ecrire
\begin{eqnarray*}
 \beta_{k}^{(\alpha)} + \frac{1}{N} \sum_{u=0}^k F_{\alpha}(\frac{u}{N})\beta_{k-u}^{(\alpha)}
&=&  \beta_{k}^{(\alpha)} \left( 1-o(1) \right) \\
&=&  \beta_{k}^{(\alpha)} \left( 1- \frac{k}{N}\right)^\alpha \left( 1+o(1) \right). 
\end{eqnarray*}
Si $k_{0}<k<N\delta _{1}$ nous pouvons consid\'erer, toujours avec le lemme \ref{lemmemoinsun}, les \'egalit\'es
\begin{eqnarray*}
  \beta_{k}^{(\alpha)} - \frac{1}{N} \sum_{u=0}^k F_{\alpha}(\frac{u}{N})\beta_{k-u}^{(\alpha)}
  &=&  \beta_{k}^{(\alpha)} - \frac{1}{N} \sum_{u=0}^k 
 \left( F_{\alpha}(\frac{u}{N})-F_{\alpha}(0) \right) \beta_{k-u}^{(\alpha)}\\
&-&F_{\alpha}(0)\frac{1}{N} \sum_{u=0}^k \beta_{k-u}^{(\alpha)}\\
&=& \beta_{k}^{(\alpha)} - \beta_{k}^{(\alpha+1)}\frac{\alpha^2}{N} 
+ \frac{1}{N^2} \sum_{u=0}^k u\beta_{k-u}^{(\alpha)}F'_{\alpha} (c_{u})
\end{eqnarray*}
avec, pour tout $u$, $ 0\le c_{u}\le \frac{u}{N}$.
Alors si $k_{0}$ suffisamment grand 
pour que l'approximation de  $\beta_{k}^{(\alpha)}$ et de  $\beta_{k}^{(\alpha+1)}$ soient pertinentes on a, en 
utilisant l'uniformit\'e de ces m\^emes approximations,
\begin{eqnarray*}
 \beta_{k}^{(\alpha)} - \beta_{k}^{(\alpha+1)}\frac{\alpha^2}{N} &=& 
\beta_{k}^{(\alpha)} - \alpha \frac{k}{N} \beta_{k}^{(\alpha)} \left( 1+o(1) \right) \\
&=& \beta_{k}^{(\alpha)} \left( 1- \frac{k}{N}\right)^\alpha \left( 1+o(1) \right) 
\end{eqnarray*}
uniform\'ement par rapport \`a $k\in [k_{0},N\delta _{1}]$.
 Si $\alpha$ positif on a
$ \frac{1}{N^2} \sum_{u=0}^k u\beta_{k-u}^{(\alpha)}F'_{\alpha} (c_{u}) =O\left( 
 \frac{1}{N^2} \sum_{u=0}^k  u \beta_{k-u}^{(\alpha)}\right),$
 et nous avons \'egalement
 \begin{eqnarray*} 
\Bigl \vert  \frac{1}{N^2} \sum_{u=0}^k u \beta_{k-u}^{(\alpha)} \Bigr \vert 
&\le & \frac{1}{N^2} \left( \sum_{u=0}^k (k-u) \vert\beta_{k-u}^{(\alpha)}\vert 
+ \sum_{u=0}^k k \vert\beta_{k-u}^{(\alpha)}\vert\right)\\
&= & O\left( \frac{k^{\alpha+1}}{N^2} \right) = O\left( k^{\alpha-1} \delta _{1}^2 \right)= o(\beta_{k}^{(\alpha)})
\end{eqnarray*}
D'autre part si $\alpha$ est n\'egatif nous pouvons \'ecrire, toujours si $k_{0}$ suffisamment grand 
pour que l'approximation des coefficients $\beta_{k}^{(\alpha)}$ soit pertinente 
$$
 \sum_{u=0}^k u \beta_{k-u}^{(\alpha)} = 
 \sum_{u=0}^{k-k_{0}} u \beta_{k-u}^{(\alpha)}+ 
 \sum_{u=k-k_{0}+1}^k u \beta_{k-u}^{(\alpha)}.
 $$
 Nous avons encore 
 $$ \frac{1}{N^2} \sum_{u=0}^{k-k_{0}} u \beta_{k-u}^{(\alpha)} 
=  \frac{1}{N^2} \left( \sum_{u=0}^{k-k_{0}} (k-u) \beta_{k-u}^{(\alpha)} 
+ \sum_{u=0}^{k-k_{0}} k \beta_{k-u}^{(\alpha)}\right).
 $$
 On a \'evidemment 
 \begin {eqnarray*}
\frac{1}{N^2}\sum_{u=0}^k (k-u) \vert\beta_{k-u}^{(\alpha)}\vert 
 &=& O( \frac{k^{\alpha+1}}{N^2})\\
 &=& O\left ( k^{\alpha-1} \delta _{1}^2 \right)= o(\beta_{k}^{(\alpha)})
\end{eqnarray*}
et 
$$ \sum_{u=0}^{k-k_{0}} k \beta_{k-u}^{(\alpha)} = -k \sum_{u=k-k_{0}+1}^\infty \beta_{k-u}^{(\alpha)} =O(\beta^{(\alpha)}_{k}).$$
Enfin en utilisant la formule d'Euler et Mac-Laurin (en supposant $k_{0}$ assez grand pour que cela ait un sens) il vient
\begin {eqnarray*}
\frac{1}{N^2} \sum_{u=k-k_{0}+1}^k  \vert\beta_{k-u}^{(\alpha)}\vert &=&
 \frac{1}{N^2} O \left( k^{\alpha+1} -  (k-k_{0})^{\alpha+1} \right)\\
 &=&  \frac{k^{\alpha+1}}{N^2} O \left( 1 -  (1-\frac{k_{0}}{k})^{\alpha+1} \right)
= O(\frac{k^{\alpha}}{N})      =o( k^{\alpha-1})
\end{eqnarray*}
ce qui ach\`eve de d\'emontrer le lemme, et l'uniformit\'e, pour $k\in [0, \delta _{0}].$\\
Reste \`a obtenir la formule pour $ k \in [N \delta _{2}, N-n_{1}]$.\\
Pour ce faire nous allons utiliser les polyn\^{o}mes orthogonaux 
$Q_{j}$, $0 \le j \le N$  associ\'es au poids $f$ et reli\'es aux polyn\^{o}mes pr\'edicteurs par 
$Q_{N} (\chi) = \overline{ P_{N}(\chi)}$. 

Nous allons utiliser \'egalement la relation  (voir \cite{I1},\cite{I2}) que 
  \begin{equation} \label{INOUE}
   \beta_{N+1,N+1} \sim  \frac{\alpha}{N}. 
   \end{equation}
D'autre part la relation (voir, par exemple \cite{SZEG})
$$ P_{m} (\chi) = P_{m-1}(\chi)  + \chi \beta_{m,m} Q_{m-1}$$
permet, en identifiant les coefficients en $N-m$ 
 $$ \beta_{N-m} = \beta_{N-m,N-1} + \beta_{N,N} \beta_{m,N-1}$$
 et pour $m \le k$ et $k$ suffisamment petit nous obtenons 
 $$ \beta_{N-k,N-m} = \beta_{N-m-(k-m),N-m} = \beta_{N-k,N-m-1}+ \beta_{N-m,N-m} \beta_{k-m,N-m-1}$$
 ce qui donne, en ajoutant cette relation pour $m$ de $0$ \`a $k$ : 
 $$ \beta_{N-k,N} = \sum_{m=0}^k \beta_{N-m,N-m} \beta_{k-m,N-m-1}.$$
 Ce qui se traduit, en utilisant le lemme (\ref{lemmemoinsun}) et l'\'equation (\ref{INOUE})
$$
\beta_{N-k,N} =(1-\frac{\alpha^2}{N}) \frac{\alpha}{N} \sum_{m=0}^k \beta_{k-m,N-m-1} + o(\frac{1}{N})
$$
 et avec les r\'esultats \'etablis au d\'ebut de la d\'emonstration 
 \begin{eqnarray*}
 \beta_{N-k,N} &=&(1-\frac{\alpha^2}{N}) \frac{\alpha}{N} \sum_{m=0}^k \left( \beta_{k-m}^{(\alpha)} - \frac{\alpha^2}{N}
 \beta _{k-m}^{(\alpha+1)}\right)+ o(\frac{1}{N})\\
 &=& (1-\frac{\alpha^2}{N})\frac{\alpha}{N} \left( \beta_{k}^{(\alpha+1)} - \frac{\alpha^2}{N} \beta_{k}^{(\alpha+2)}\right)
 + o(\frac{1}{N})
 \end{eqnarray*}
 Ce qui donne 
 \begin{equation} \label{INOUE2}
 \beta_{N-k+1,N} = \frac{\alpha \beta_{k}^{(\alpha+1)} }{N} +  o(\frac{1}{N}).
 \end{equation}
 Ce qui s'\'ecrit aussi si nous consid\'erons un entier $j \in [N\delta _{2},N-n_{1}]$
 \begin{eqnarray*}
g_{1}(1)\beta_{j+1,N} &=& \frac{\alpha (N-j) ^\alpha }{N \Gamma (\alpha+1)} = \frac{1}{\Gamma(\alpha)} ( 1-\frac{j}{N} )^\alpha
N ^{\alpha-1}+o(\frac{1}{N})\\
&=& \frac{1}{\Gamma(\alpha)} ( 1-\frac{j}{N} )^\alpha( \frac{N-j+j} {j} ) ^{\alpha-1} j^{\alpha-1} +o(\frac{1}{N})\\
 \end{eqnarray*}
 Soit 
 $$\beta_{j+1,N} = \beta_{j}^{(\alpha)} ( 1-\frac{j}{N} )^\alpha +o(\frac{1}{N}).$$

 L'uniformit\'e du r\'esultat est alors assur\'ee par la fa\c con dont $\delta _{2}$ tend vers $1$.
 \section {D\'emonstration du th\'eor\`eme principal}
On doit calculer la trace de  $ \left( T_{N} (f_{1}) T_{N}^{-1} (f_{2}) \right)$ avec 
$f_{1}= \vert 1- \chi \vert ^{2 \alpha_{1}} c_{1}$ et $f_{2}= \vert 1- \chi \vert ^{2 \alpha_{2}} c_{2}$ 
o\`u $-\frac{1}{2} < \alpha_{1}, \alpha_{2} < \frac{1}{2}$ et o\`u 
$ c_{1}$ et $c_{2}$ sont des fonctions r\'eguli\`eres.\\
On peut \'ecrire 
$$ \Tr  \left( T_{N} (f_{1}) T_{N}^{-1} (f_{2}) \right) = 
\hat f_{1} (0) \sum _{k=0}^N \left(T_{N}^{-1} (f_{2})\right)_{k,k}
+ 2 \Re \left( \sum_{s=1}^N \hat f_{1} (s) \sum_{k=0}^{N-s} T_{N}^{-1} (f_{2})_{k,k+s} \right),$$
et aussi (avec $s>0$) et en utilisant le lemme \ref{Gohberg}
$$  \sum_{k=0}^{N-s} T_{N}^{-1} (f_{2})_{k,k+s} = (N-1-s) A_{1} (s) +A_{2} (s) $$
avec $$ \left \lbrace
\begin{array} {l}
A_{1} (s) = \sum_{l=0}^{N-s} \beta_{l,N} \overline{\beta_{l+s,N} } \\
A_{2} (s) = -2  \left(\sum_{l=0}^{N-s} l \beta_{l,N} \overline{\beta_{l+s,N} } \right)
\end{array}
\right.$$
o\`u les $\beta_{l,N}$, $0 \le l \le N$, sont les coefficients du polyn\^{o}me pr\'edicteur de la fonction $f_{2}$.
Nous avons \'etabli plus haut que si $n_{0}$ est un entier tel que $\frac{n_{0}} {N} =o(1)$ nous pouvons 
\'ecrire 
\begin{equation} \label{rem}
 \beta_{l,N} = \beta_{l}^{(\alpha_{2})} (1- \frac{l}{N} )^{\alpha_{2}} \left (1+o(1) \right)
 \end{equation}
  uniform\'ement pour 
tout entier $l$ dans $[0, N-n_{0}]$. 
Dans un premier temps nous allons d\'emontrer le lemme suivant 
\begin{lemme}\label{lemme1}
Avec les hypoth\`eses du th\'eor\`eme nous avons quel que soit $s$, $0\le s\le N-n_0$, uniform\'ement par rapport \`a $s$ 
\begin{enumerate}
\item
Si $\alpha_{2} \in ]-\frac{1}{2},0[$ alors 
\begin{eqnarray*}
\sum _{l=0}^{N-s} T_{N}^{-1} (f_{2}) _{l,l+s} &=& N (1-\frac{s}{N} )^{\alpha_{2}+1}
 \widehat {\frac{1}{f_{2}}} (-s) 
+ (1-\frac{s}{N} )^{\alpha_{2}} \tau_2 \left( \frac{s}{N} \right) \left( \sum _{u\ge 0} u \beta_{u}^{\alpha_{2}} 
\overline{\beta_{u}^{\alpha_{2}}}\right)
 +\\ 
&+&
\frac{ N^{2\alpha_{2}} }{\Gamma^2 (\alpha_{2})c_{2}(1)}   (1- \frac{s}{N} ) ^{2\alpha_{2}} 
  F_{1}(\frac{s}{N}) + o(N^{2\alpha_2})
\end{eqnarray*}
o\`u $F_{1}$ est une fonction continue sur $[0,1]$,
et $\tau_2 (t) = -\left(\alpha_2+\left((1-t) (\alpha_2+2)\right)\right).$\\
\item
Si $\alpha_{2}\in ]0,  \frac{1}{2}[$ alors 
$$ \sum _{l=0}^{N-s} T_{N}^{-1} (f_{2}) _{l,l+s} = N (1-\frac{s}{N} )^{\alpha_{2}+1}
 \widehat {\frac{1}{f_{2}}} (-s) + \frac{ N^{2\alpha_{2}} }{\Gamma^2 (\alpha_{2})c_{2}(1)} F_2 (\frac{s}{N}) +o(N^{2\alpha_{2}})$$
 o\`u $F_2$ est une fonction continue sur $[0,1]$,
   \end{enumerate}
\end{lemme}

\subsection{ D\'emonstration du lemme \ref{lemme1}}
\begin{remarque}
La d\'emonstration de l'uniformit\'e des restes, qui est un peu fastidieuse, a \'et\'e repouss\'e dans l'appendice.
\end{remarque}
\subsubsection{Calcul du coefficient $A_{1}$}
En supposant que $n_{1}$ est un entier tel que $\frac{n_{1}} {N} =o(1)$ et $n_{1} \le n_{0}$ nous pouvons \'ecrire la d\'ecomposition :
$$ A_{1}(s)= \sum_{l=0}^{N-s-n_{1}} \beta_{l,N} \overline{\beta_{l+s,N} } + 
\sum_{l=N-s-n_{1}+1}^{N-s} \beta_{l,N} \overline{\beta_{l+s,N} }.$$ 
Compte tenu de la remarque \ref{rem} nous pouvons \'ecrire, en posant $s' = s+n_{1}$  
$$
 \sum_{l=0}^{N-s'} \beta_{l,N} \overline{\beta_{l+s,N} } 
 =\sum_{l=0}^{N-s'}  \beta_{l}^{(\alpha_{2})} \overline{\beta_{l+s} ^{(\alpha_{2})} }
(1- \frac{l}{N} )^{\alpha_{2}} \overline {(1- \frac{l+s}{N} )^{\alpha_{2}}} 
\left(1+o(1)\right).
$$
Nous avons alors
\begin{eqnarray*}
 & & \sum_{l=0}^{N-s'}  \beta_{l}^{(\alpha_{2})} \overline{\beta_{l+s} ^{(\alpha_{2})} }
(1- \frac{l}{N} )^{\alpha_{2}} \overline {(1- \frac{l+s}{N} )^{\alpha_{2}}} =\\
&=& 
 \sum_{l=0}^{N-s'}  \beta_{l}^{(\alpha_{2})} \overline{\beta_{l+s} ^{(\alpha_{2})} }
 \left(\left( (1 - \frac{l}{N} ) ^{\alpha_{2}} -1 \right) +1\right) 
 \left( \left( (1 - \frac{l+s}{N} ) ^{\alpha_{2}}-(1 - \frac{s}{N} ) ^{\alpha_{2}}\right) + (1 - \frac{s}{N} ) ^{\alpha_{2}}\right)
 \end{eqnarray*}
En d\'eveloppant on a 
$$ 
\sum_{l=0}^{N-s'}  \beta_{l}^{(\alpha_{2})} \overline{\beta_{l+s} ^{(\alpha_{2})} }
(1- \frac{l}{N} )^{\alpha_{2}} \overline {(1- \frac{l+s}{N} )^{\alpha_{2}}} 
=
A'_{1} + A'_{2}+ A'_{3} + A'_{4},
$$
avec
$$
\left \lbrace 
\begin{array}{l} 
A'_{1}= \displaystyle{ \left( \sum_{l=0}^{N-s'}  \beta_{l}^{(\alpha_{2})} \overline{\beta_{l+s} ^{(\alpha_{2})} }\right)
(1 - \frac{s}{N} ) ^{\alpha_{2}}}\\
A'_{2} = \displaystyle{\left( \sum_{l=0}^{N-s'}  \beta_{l}^{(\alpha_{2})} \overline{\beta_{l+s} ^{(\alpha_{2})} }
\left((1 - \frac{l}{N} ) ^{\alpha_{2}} -1\right)\right) (1 - \frac{s}{N} ) ^{\alpha_{2}}}\\
A'_{3} =\displaystyle{ \sum_{l=0}^{N-s'}  \beta_{l}^{(\alpha_{2})} \overline{\beta_{l+s} ^{(\alpha_{2})} }
 \left( (1 - \frac{l+s}{N} ) ^{\alpha_{2}}-(1 - \frac{s}{N} ) ^{\alpha_{2}}\right)} \\
 A'_{4}=  \displaystyle{\sum_{l=0}^{N-s'}  \beta_{l}^{(\alpha_{2})} \overline{\beta_{l+s} ^{(\alpha_{2})} }
\left( (1 - \frac{l}{N} ) ^{\alpha_{2}} -1 \right) 
 \left( (1 - \frac{l+s}{N} ) ^{\alpha_{2}}-(1 - \frac{s}{N} ) ^{\alpha_{2}}\right)}.
\end{array}
\right.
$$
Nous avons 
\begin{eqnarray*}
A'_{1} &=& (1 - \frac{s}{N} ) ^{\alpha_{2}}
\left( \widehat {\frac{1}{f_{2}} }(-s) - \sum _{N-s'+1} ^{+ \infty} \beta_{l}^{(\alpha_{2})} \overline{\beta_{l+s} ^{(\alpha_{2})} } \right)\\
&=& (1 - \frac{s}{N} ) ^{\alpha_{2}}
\left( \widehat {\frac{1}{f_{2}} }(-s) -\right.\\
&-&\left.\frac{N^{2\alpha_{2}-1}}{c_{2}(1) \Gamma ^2 (\alpha_{2})}
\int _{1-s'/N} ^{+ \infty} t^{\alpha_{2} -1} (t+\frac{s}{N} )^{\alpha_{2} -1} dt \right)
+o(N^{2\alpha_2-1}).
\end{eqnarray*}
On peut remarquer que quand $\frac{s}{N}$ tend vers $1$ et $\alpha_{2}$ n\'egatif nous pouvons \'ecrire 
$$ A'_{1}  \sim (1- \frac{s}{N}) ^{\alpha_{2}}  \widehat {\frac{1}{f_{2}} }(-s) 
-  (1- \frac{s}{N}) ^{2 \alpha_{2}} \frac{N^{2\alpha_{2}}}{ \alpha_{2} c_{2}(1) \Gamma ^2 (\alpha_{2})}+ o(N^{2\alpha_{2}})
$$
Nous avons d'autre part \\
\begin{itemize}
\item [i)]
Si $\alpha_{2} \in ]- \frac{1}{2},0[$
\begin{eqnarray*}
A'_{2} &=&
 \left( \sum_{l=0}^{N-s'}  \beta_{l}^{(\alpha_{2})} \overline{\beta_{l+s} ^{(\alpha_{2})} }
\left( (1 - \frac{l}{N} ) ^{\alpha_{2}} -1 + \alpha_{2} \frac{l}{N} \right)\right) (1- \frac{s}{N}) ^{ \alpha_{2}} - \\
&-& \frac{\alpha_{2} }{N}
\left(\sum_{l=0}^{N-s'}  \beta_{l}^{(\alpha_{2})} \overline{\beta_{l+s} ^{(\alpha_{2})} } l  \right)
(1- \frac{s}{N}) ^{ \alpha_{2}}.
\end{eqnarray*}
Ce qui donne 
\begin{eqnarray*}
&&\frac{\alpha_{2} }{N}
\left(\sum_{l=0}^{N-s'}  \beta_{l}^{(\alpha_{2})} \overline{\beta_{l+s} ^{(\alpha_{2})} } l  \right)
= 
\frac{\alpha_{2} }{N} 
 \left( \sum _{u\ge 0} u \beta_{u}^{(\alpha_{2})} 
\overline{\beta_{u+s}^{(\alpha_{2})}}\right)
 -\\
&-&  N^{2\alpha_{2} -1} 
\frac{\alpha_{2} }{ c_{2}(1) \Gamma ^2(\alpha_{2})} 
\int_{1-s'/N}^{+ \infty} t^{\alpha_{2}} (t+\frac{s}{N} )^{\alpha_{2}-1} dt 
+o(N^{2\alpha_2-1}),
\end{eqnarray*}
et d'autre part 
$$
\left( \sum_{l=0}^{N-s'}  \beta_{l}^{(\alpha_{2})} \overline{\beta_{l+s} ^{(\alpha_{2})} }
\left( (1 - \frac{l}{N} ) ^{\alpha_{2}} -1 + \alpha_{2} \frac{l}{N} \right)\right)
=
 \frac{N^{2\alpha_{2} -1}  }{ c_{2}(1) \Gamma ^2(\alpha_{2})} \Phi_{1}( \frac{s'}{N}) + o(N^{\alpha_2-1})$$
 avec 
 $$ \Phi_{1}( \frac{s'}{N}) =
\int _{0}^{1-s'/N} t^{\alpha_{2}-1} (t+ \frac{s}{N})^{\alpha_{2}-1} \left( (1-t) ^{\alpha_{2}}-1+ \alpha_{2} t \right) dt.
$$
En posant finalement 
$$ \Phi_{2}( \frac{s'}{N}) = 
\Phi_{1}( \frac{s'}{N})  + \alpha_2 \int_{1-s'/N}^{+ \infty} t^{\alpha_{2}} (t+\frac{s}{N} )^{\alpha_{2}-1} dt
$$ 
nous pouvons \'ecrire 
\begin{align*}
 A'_{2} &= \left(
-\frac{\alpha_{2} }{N} 
 \left( \sum _{u\ge 0} u \beta_{u}^{(\alpha_{2})} 
\overline{\beta_{u+s}^{(\alpha_{2})}}\right)
+\right.\\
&\left.+ N^{2\alpha_{2}-1} \frac{1}{ c_{2}(1) \Gamma ^2(\alpha_{2})} \Phi_{2}( \frac{s'}{N})\right) (1-\frac{s}{N} )^{\alpha_{2}}+ o(N^{\alpha_2-1}),
\end{align*}
o\`u $\Phi_{2}$ est une fonction continue sur $[0,1]$.
\item [ii)]
Si $\alpha_{2} \in ]0, \frac{1}{2}[$ nous avons 
$$ A'_{2} = (1-\frac{s}{N} )^{\alpha_{2}} N^{2\alpha_{2}-1} \frac{1}{ c_{2}(1) \Gamma ^2(\alpha_{2})} 
\int_{0}^1 t^{\alpha_{2}} (t+\frac{s}{N} ) ^{\alpha_{2}-1} \left( (1-t)^{\alpha_{2}} -1 \right) dt +o(N^{\alpha_2 -1})$$
\end{itemize}
Calculons maintenant le terme $A'_{3}$.L\`a aussi nous devons distinguer les cas $\alpha_{2}$ positif et $\alpha_{2}$ n\'egatif. 
\begin{itemize}
\item [i)]
Si $\alpha_{2 } \in ]-\frac{1}{2}, 0[$.\\
Nous pouvons \'ecrire  
\begin {eqnarray*}
A'_{3} &=& \sum_{l=0} ^{N-s'} \beta_{l}^{(\alpha_{2})} \overline{\beta_{l+s}^{(\alpha_{2})}}
\left( (1-\frac{l+s}{N})^{\alpha_{2}} -(1-\frac{s}{N}) ^{\alpha_{2}} \right) \\
&= & \sum_{l=0} ^{N-s'} \beta_{l}^{(\alpha_{2})} \overline{\beta_{l+s}^{(\alpha_{2})}}
\left( (1-\frac{l+s}{N})^{\alpha_{2}} -(1-\frac{s}{N}) ^{\alpha_{2} }+ \alpha_{2} \frac{l}{N} (1-\frac{s}{N}) ^{\alpha_{2}-1}
\right) -\\
&-&\frac{\alpha_{2}}{N} \sum_{l=0} ^{N-s'}  l \beta_{l}^{(\alpha_{2})} \overline{\beta_{l+s}^{(\alpha_{2})}}
 (1-\frac{s}{N} )^{\alpha_{2}-1} .
\end {eqnarray*}
Nous avons 
\begin {eqnarray*}
& & \frac{\alpha_{2}}{N}\sum_{l=0} ^{N-s'}l  \beta_{l}^{(\alpha_{2})} \overline{\beta_{l+s}^{(\alpha_{2})}}
 (1-\frac{s}{N}) ^{\alpha_{2} -1} =  (1-\frac{s}{N}) ^{\alpha_{2} -1} \frac{\alpha_{2}}{N}
\sum_{u=0}^\infty u \beta_{u}^{(\alpha_{2})} \overline{\beta_{u+s}^{(\alpha_{2})}} -\\
&-& (1-\frac{s}{N}) ^{\alpha_{2} -1}
 \frac{\alpha_{2}N^{2\alpha_{2}-1}} { c_{2}(1) \Gamma^2(\alpha_{2})}\int _{1-s'/N}^{+ \infty} 
t^{\alpha_{2}} (t+\frac{s}{N} )^{\alpha_{2}-1} dt+o(N^{2\alpha_2-1})
\end {eqnarray*}
et de m\^{e}me 
\begin{eqnarray*}
& &\sum_{l=0} ^{N-s'}   \beta_{l}^{(\alpha_{2})} \overline{\beta_{l+s}^{(\alpha_{2})}}
\left( (1-\frac{l+s}{N})^{\alpha_{2} } -  (1-\frac{s}{N}) ^{\alpha_{2}} + \alpha_{2} \frac{l}{N} (1-\frac{s}{N}) 
^{\alpha_{2}-1} \right) = \\
 &=& \frac{N^{2\alpha_{2}-1}} { c_{2}(1) \Gamma^2(\alpha_{2})} \Phi_{3}(\frac{s}{N}) +o(N^{2\alpha_2-1}),
 \end{eqnarray*}
 Avec 
 $$ \Phi_{3}(\frac{s}{N}) = \int _{0} ^{1-s'/N} t^{\alpha_{2}-1} (t+ \frac{s}{N} )^{\alpha_{2}-1}
 \left( (1-t-\frac{s}{N}) ^{\alpha_{2}} - (1-\frac{s}{N} )^{\alpha_{2}} 
 + \alpha_{2} t (1-\frac{s}{N} )^{\alpha_{2}-1}\right) dt.$$
En \'ecrivant 
 \begin{eqnarray*}
& &  t^{\alpha_{2}-1} (t+ \frac{s}{N} )^{\alpha_{2}-1}
 \left( (1-t-\frac{s}{N}) ^{\alpha_{2}} - (1-\frac{s}{N} )^{\alpha_{2}} 
 + \alpha_{2} t (1-\frac{s}{N} )^{\alpha_{2}-1}\right)\\
 &=& (t+ \frac{s}{N} )^{\alpha_{2}-1} \displaystyle{ \sum_{n\ge 2} \delta_{n} t^{\alpha_{2}+n -1} 
 (1- \frac{s}{N} )^{\alpha_{2}-n}}
 \end{eqnarray*}
 o\`u les coefficients $\delta _{n}$ sont les coefficients du d\'eveloppement  enti\`ere de la fonction 
 $ t \rightarrow (1-t)^{\alpha_{2}}$ 
 nous pouvons conclure que 
 $\Phi_{3} (u) \sim O \left(1- u) ^{2 \alpha_{2}}\right)$ au voisinage de $1$.
 En posant 
 $$ \Phi_{4} (\frac{s }{N} ) = \alpha_{2} \int_{1-s'/N}^{ +\infty} t^{2\alpha_{2}} (t+ \frac{s}{N}) ^{\alpha_{2}-1} dt$$
 nous pouvons \'ecrire finalement 
 \begin{eqnarray*} 
 A'_{3}&=& (1-\frac{s}{N} )^{\alpha_{2}-1} \left( -\frac{\alpha_2}{N}  \left( \sum _{u\ge 0} u \beta_{u}^{(\alpha_{2})} 
\overline{\beta_{u+s}^{(\alpha_{2})}}\right)
\right)+\\
 &+&
 \frac{N^{2\alpha_{2}-1}} { c_{2}(1) \Gamma^2(\alpha_{2})}  
   \left(  (1- \frac{s}{N} ) ^{\alpha_{2}-1 }\Phi_{4} (\frac{s}{N})+\Phi_{3}\left( \frac{s}{N} \right)\right)+o(N^{2\alpha_2-1})
 \end{eqnarray*}
 o\`u $\Phi_{3}$ est une fonction continue sur $[0,1[$ et \'equivalente \`a 
 $O \left( (1- t) ^{2 \alpha_{2}}\right)$ en $1$, et o\`u $\Phi_{4} $ est une fonction continue sur $[0,1]$.
 \item [ii)]
 Si $\alpha_{2}\in ]0, \frac{1}{2}[$ nous pouvons \'ecrire de m\^{e}me
 $$ \frac{\alpha_{2}}{N}\sum_{l=0} ^{N-s'}l  \beta_{l}^{(\alpha_{2})} \overline{\beta_{l+s}^{(\alpha_{2})}}
=
 \frac{\alpha_{2}N^{2\alpha_{2}-1}} { c_{2}(1) \Gamma^2(\alpha_{2})}\int _0^1
t^{\alpha_{2}} (t+\frac{s}{N} )^{\alpha_{2}-1} dt +o(N^{2\alpha_2-1})
$$
et, comme ci-dessus 
\begin{eqnarray*}
& &\sum_{l=0} ^{N-s'} l  \beta_{l}^{(\alpha_{2})} \overline{\beta_{l+s}^{(\alpha_{2})}}
\left( (1-\frac{l+s}{N})^{\alpha_{2} } -  (1-\frac{s}{N}) ^{\alpha_{2}} + \alpha_{2} \frac{l}{N} (1-\frac{s}{N}) 
^{\alpha_{2}-1} \right) = \\
 &=& \frac{N^{2\alpha_{2}-1}} { c_{2}(1) \Gamma^2(\alpha_{2})} \Phi_{3}(\frac{s'}{N})+o(N^{2\alpha_2-1}). 
 \end{eqnarray*}
D'o\`u 
$$
 A'_{3}=  \left(  \frac{N^{2\alpha_{2}-1}} { c_{2}(1) \Gamma^2(\alpha_{2})}  
   \left( (1-\frac{s}{N} )^{\alpha_{2}-1}\tilde \Phi_{4} (\frac{s'}{N})+\Phi_{3}\left( \frac{s'}{N} \right)\right)\right)+o(N^{2\alpha_2-1})
,
 $$
avec 
$$\tilde \Phi_{4}(\frac{s}{N}) = \alpha_{2} \int_0^1t^{2\alpha_{2}} (t+ \frac{s}{N}) ^{\alpha_{2}-1} dt$$
 \end{itemize}
 Il est d'autre part facile de se convaincre que 
\begin{eqnarray*}
 A'_{4} &=& \frac{N^{2\alpha_{2}-1}} { c_{2}(1) \Gamma ^2 (\alpha_{2}) }
 \int _{0} ^{1-\frac{s'}{N} } t^{\alpha_{2}-1} \left(( t+ \frac{s}{N})^{\alpha_{2}-1} \right) \\
& &  \left( (1- t-\frac{s}{N} )^{\alpha_{2}} -(1-\frac{s}{N})^{\alpha_{2}}\right) \left( (1-t)^{\alpha_{2}-1} \right)+o(N^{2\alpha_2-1})
 \\
 &=&  \frac{N^{2\alpha_{2}-1}} { c_{2}(1) \Gamma ^2 (\alpha_{2})} \Phi_{5}(\frac{s}{N}) +o(N^{2\alpha_2-1})
 \end{eqnarray*} 
o\`u $ \Phi_{5}$ est une fonction continue sur $[0,1[$ et o\`u $\Phi_{5} (u)$ est \'equivalent \`a 
  $O \left( (1- u) ^{2 \alpha_{2}+1}\right)$ en $1$. \\
   Enfin en utilisant l'\'equation \ref{INOUE2}
      qui est apparue dans la d\'emonstration du th\'eor\`eme \ref{theoremezero} on obtient que, quel que soit le signe de $\alpha_{2}$
      on a   
        $$\sum_{l=N-s-n_{1}}^{N-s} \overline{ \beta_{l,N} \beta _{l+s,N}}= O(N^{\alpha_{2}-2})$$
        ce qui implique 
        $$ (N-s-1) \sum_{l=N-s-n_{1}}^{N-s} \overline{ \beta_{l,N} \beta _{l+s,N}} = o (N^{2\alpha_{2}}).$$
  En regroupant on peut donc \'ecrire 
  \begin{itemize}
  \item[$\bullet$]
  Si $ \alpha_{2} \in ]- \frac{1}{2}, 0[$
  \begin{eqnarray*} 
 ( N+1-s) A_{1} (s)& = & N (1-\frac{s}{N} )^{\alpha_2+1}\widehat{ \frac {1} {f_{2}} }(-s) 
 - \alpha_{2 } (2 - \frac{s}{N})  (1- \frac{s}{N} )^{\alpha_{2}}
 \left( \frac{s}{N} \right) \left( \sum _{u\ge 0} u \beta_{u}^{(\alpha_{2})} 
\overline{\beta_{u+s}^{(\alpha_{2})}}\right)
 +\\
  &+& \frac{N^{2\alpha_{2}} }{ c_{2}(1) \Gamma ^2 (\alpha_{2})} 
  (	1-\frac{s}{N} )^{2\alpha_{2}} \tilde F_{1}(\frac{s}{N}) +o(N^{2\alpha_2})
  \end{eqnarray*}
  avec $\tilde F_{1}$ une fonction continue sur $[0,1]$. 
  \item[$\bullet \bullet$]
  Si $ \alpha_{2} \in ]0, \frac{1}{2}[$
  $$ ( N+1-s) A_{1} (s)  = N (1-\frac{s}{N} )^{\alpha_2+1}\widehat{ \frac {1} {f_{2}} }(-s) 
  +\frac{N^{2\alpha_{2}} }{ c_{2}(1) \Gamma ^2 (\alpha_{2})} 
\tilde F_{2}(\frac{s}{N})+o(N^{2\alpha_2}),$$
o\`u $\tilde F_{2}$ est une fonction continue sur $[0,1]$.
\end{itemize}
    \subsubsection{ Calcul du coefficient $A_{2}$}
    Toujours avec la m\^{e}me d\'efinition de $s'$ on a 
    Nous avons $A_{2}(s)= -2 \sum _{l=0}^{N-s'} l \beta_{l,N} \overline{\beta_{l+s,N}}.$
    \begin{itemize}
    \item[i)]
    Dans le cas o\`u $\alpha_{2} \in ]- \frac{1}{2},0[ $ nous avons  
 \begin{eqnarray*}
 A_{2}(s) &=&-2  \sum _{l=0}^{N-s'} l \beta_{l,N} \overline{\beta_{l+s,N}} \\
 &=&-2\sum _{l=0}^{N-s'}  l \beta_{l}^{(\alpha_{2})} \overline{\beta_{l+s}^{\alpha_{2}}} (1-\frac{l}{N})^{\alpha_{2}} 
 (1-\frac{l+s}{N})^{\alpha_{2}} \\
 &=& -2\sum _{l=0}^{N-s'}  l \beta_{l}^{(\alpha_{2})} \overline{\beta_{l+s}^{\alpha_{2}}} 
 \left( (1-\frac{l}{N})^{\alpha_{2}} (1-\frac{l+s}{N})^{\alpha_{2}} - (1-\frac{s}{N})^{\alpha_{2}}\right) -\\
 &-&2 (1-\frac{s}{N})^{\alpha_{2}} \sum _{l=0}^{N-s'}  l \beta_{l}^{(\alpha_{2})} \overline{\beta_{l+s}^{(\alpha_{2})}}. 
 \end{eqnarray*}
Soit, finalement :
  \begin{eqnarray*}
 A_{2}(s) &=&-2 (1-\frac{s}{N})^{\alpha_{2}}  \left( \sum _{u\ge 0} u \beta_{u}^{(\alpha_{2})} 
\overline{\beta_{u+s}^{(\alpha_{2})}}\right)
 + \\
 &+& 2 (1-\frac{s}{N})^{\alpha_{2}}  \frac{N^{2\alpha_{2}} }{ c_{2}(1) \Gamma ^2 (\alpha_{2})} 
\int _{1-s'/N} ^{+ \infty} t^{\alpha_{2}} (t+ \frac{s}{N} ) ^{\alpha_{2}-1} dt \\
&-& 2  \frac{N^{2\alpha_{2}} }{ c_{2}(1) \Gamma ^2 (\alpha_{2})} \int_{0}^{1-s'/N} t^{\alpha_{2}} (t+ \frac{s}{N} ) ^{\alpha_{2}-1} 
\left( (1-t)^{\alpha_{2}} (1-t-\frac{s}{N})^{\alpha_{2}} -(1-\frac{s}{N})^{\alpha_{2}}  \right) dt\\
&+& o(N^{2\alpha_2}). 
 \end{eqnarray*}
 En posant 
  \begin{eqnarray*}
  H_{3} (\frac{s}{N}) &= &-\int _{1-s'/N} ^{+ \infty} t^{\alpha_{2}} (t+ \frac{s}{N} ) ^{\alpha_{2}}dt +\\
 H_{4} (\frac{s}{N}) &= &\int_{0}^{1-s'/N} t^{\alpha_{2}} (t+ \frac{s}{N} ) ^{\alpha_{2}} 
\left( (1-t)^{\alpha_{2}} (1-t-\frac{s}{N})^{\alpha_{2}} -(1-\frac{s}{N})^{\alpha_{2}}  \right) dt,
\end{eqnarray*}
  nous pouvons \'ecrire 
  \begin{align*}
   A_{2}(s) &= -2(1-\frac{s}{N})^{\alpha_{2}} \left( \left( \sum _{u\ge 0} u \beta_{u}^{(\alpha_{2})} 
\overline{\beta_{u+s}^{(\alpha_{2})}}\right)
 + \frac{N^{2\alpha_{2}} }{ c_{2}(1) \Gamma ^2 (\alpha_{2})}
  H_{3} (\frac{s}{N})\right)-\\ 
   &- \frac{N^{2\alpha_{2}} }{ c_{2}(1) \Gamma ^2 (\alpha_{2})} 
   H_{4} (\frac{s}{N}) +o(N^{2\alpha_2})
   \end{align*}
   
  o\`u $H_{3}$ et $H_4$ sont des fonctions continues sur $[0,1]$ 
  (on vrifie aisment que $H_4 (u) = O\left( (1-u)^{\alpha_2+1} \right)$
  au voisinage de $1$. 
  \item[ii)]
  Le cas $\alpha_{2}\in ]0, \frac{1}{2}[$ se traite plus rapidement. Il suffit de remarquer que :
  
    \begin{align*} 
    A_{2}(s) &=-2 \sum _{l=0}^{N-s'}  l \beta_{l}^{(\alpha_{2})} \overline{\beta_{l+s}^{\alpha_{2}}} (1-\frac{l}{N})^{\alpha_{2}} 
 (1-\frac{l+s}{N})^{\alpha_{2}} \\
 &= -2\frac{N^{2\alpha_{2}} }{ c_{2}(1) \Gamma ^2 (\alpha_{2})} \int_{0}^{1} 
  t^{\alpha_{2}} (t+\frac{s}{N})^{\alpha_{2}-1} (1-t)^{\alpha_{2}} (1-t - \frac{s}{N})^{\alpha_{2}} dt +o(N^{2\alpha_2}).
        \end {align*}
        \end{itemize}
       Pour r\'esumer nous pouvons poser 
           \begin{itemize}
       \item[$\bullet$] 
       Si $\alpha_{2}<0$ alors 
        \begin{align*}
  A_{2}(s) &=    -2(1-\frac{s}{N})^{\alpha_{2}}  \left( \frac{s}{N} \right) \left( \sum _{u\ge 0} u \beta_{u}^{(\alpha_{2})} 
                      \overline{\beta_{u+s}^{(\alpha_{2})}}\right)\\
           &- \frac{N^{2\alpha_{2}} }{ c_{2}(1) \Gamma ^2 (\alpha_{2})}
            \left(-2(1-\frac{s}{N})^{\alpha_{2}} H_{3} (\frac{s}{N})+ H_{4} (\frac{s}{N}) 
\right)+o(N^{2\alpha_2}). 
        \end{align*}   
         \item[$\bullet \bullet$]
       Si $\alpha_{2}>0$ alors 
       $$ A_{2}(s) = -2 \frac{N^{2\alpha_{2}} }{ c_{2}(1) \Gamma ^2 (\alpha_{2})} H_{5}(\frac{s}{N}) +o(N^{2\alpha_{2}}).$$
       \end{itemize}
  D'autre part quelque soit le signe de $\alpha_{2}$ , nous avons, comme pour le calcul de $A_{1}$ nous 
  avons 
  $$ -2\sum_{l=N-s-n_{1}}^{N-s} l \beta_{l,N} \overline{ \beta_{l+s,N}} =o(N^{2\alpha_{2}}).$$
  En r\'eunissant les r\'esultats des points abord\'es nous obtenons l'\'enonc\'e du lemme \ref{lemme1}
  \subsection {D\'emonstation du th\'eor\`eme dans le cas o\`u $\alpha_{2} \in ]- \frac{1}{2},0[$}
  Il s'agit de calculer la somme 
  $$\widehat {f_{1} }(0) \trace \left(T_{N} (f_{2})\right)^{-1} +2 \Re \left(
   \sum_{s=1} ^{N} \widehat{f_{1}} (s) \sum_{l=0}^{N-s} 
   \left(T_{N} (f_{2})\right)^{-1}_{l,l+s} \right).$$
   Rappelons que dans le cas o\`u l'exposant $\alpha_{2}$ est n\'egatif il a \'et\'e \'etabli que 
   $$ \trace \left(T_{N}^{-1} (f_{2})\right)) = (N+1) \widehat{\left(\frac{1}{f_{2}}\right)} (0) 
   + \langle \ln f_{2} \vert \frac{1}{f_{2}}\rangle_{2,1/2}  +o(1)$$
  o\`u $\langle \, \vert \rangle_{2,1/2} $ d\'esigne le produit scalaire dans $A(\mathbb {T} , \frac{1}{2}) =\{ \rho \in L^2 ( \mathbb {T}) / \displaystyle{\sum_{m\in \mathbb Z} \vert m \vert \vert \hat \rho (m) \vert ^ 2 < \infty} \} $
   (voir  \cite{RS09} ). 
   En utilisant toujours la propri\'et\'e \ref{INOUE2} extraite de la d\'emonstration du th\'eor\`eme \ref{theoremezero} 
  il vient, avec $\vert N-s\vert < n_{0}$, avec $n_{0}$ comme dans le pr\'eambule de la d\'emonstration, et si $n_0$ assez petit :
  $$ \sum_{l=0}^{N-s} \beta_{l,N} \overline {\beta_{l+s,N}} =O \left( \frac{1}{N}\right) $$
  et
 $$ \sum_{l=0}^{N-s}l \beta_{l,N} \overline {\beta_{l+s,N}} =O \left( \frac{1}{N}\right). $$
 Soit 
 $$ \sum_{s=N-n_{0}} ^N \widehat f_{1}(s) (N-s) \sum_{l=0}^{N-s} \beta_{l,N} \overline {\beta_{l+s,N}} = O(N^{-2\alpha_{2}-2})$$
 et 
 $$ \sum_{s=N-n_{0}} ^N \widehat f_{1}(s) \sum_{l=0}^{N-s}l \beta_{l,N} \overline {\beta_{l+s,N}} =O(N^{-2\alpha_{2}-2}).$$
 Nous allons donc \'etudier ne fait la somme 
 $$\widehat {f_{1} }(0) \trace \left(T_{N} (f_{2})\right)^{-1} +2 \Re \left(
   \sum_{s=1} ^{N'} \widehat{f_{1}} (s) \sum_{l=0}^{N'-s} 
   \left(T_{N} (f_{2})\right)^{-1}_{l,l+s} \right)$$ 
   avec $N'=N-n_{0}.$ Nous v\'erifierons ensuite que l'approximation trouv\'ee est d'ordre sup\'erieur \`a $O(N^{-2\alpha_{2}-2})$.\\
  D'apr\`es le lemme pr\'ec\'edent la somme pr\'ec\'edente se d\'ecompose en quatre sommes que nous 
  allons traiter s\'epar\'ement.
  \subsubsection {Calcul de $ (N+1) \left(\hat{f_{1}} (0) \widehat{\left(\frac{1}{f_{2}}\right)} (0)+ 
  \displaystyle{2 \sum _{s=1} ^{N^{\prime}} \widehat {f_1} (s) (1- \frac{s}{N}) ^ {\alpha_2 +1} \widehat{ f_{2}^{-1} } (-s)}\right) $}
  Il est connu que pour $s$ suffisamment grand on a
  $$ \widehat {f_1}(s) = \frac{{c_1} (1) } { \Gamma (-2\alpha_1)} s^{-2\alpha_1-1} 
  +o(s^{-2\alpha_1-1}) ,
  \quad  \widehat {f_{2}^{-1}}(s) = \frac{1}{c_2 (1)  \Gamma (2\alpha_2)} 
  s^{2\alpha_2-1}+ o(s^{-2\alpha_2-1}).$$
   Trois cas sont  distinguer pour calculer cette somme. \\
  a) $ \frac{1}{2} > \alpha_2 - \alpha_1 >0.$\\
  Sous cette hypothse nous pouvons \'ecrire,  en posant  
  $C_1 = \frac{c_1(1)}{c_2(1)} \frac{1}{\Gamma (-2\alpha_1)\Gamma (2 \alpha_2)}$ 
  \begin{align*} 
   (N+1) \sum _{s=1} ^N \widehat {f_1} (s) (1- \frac{s}{N}) ^
  {\alpha_2 +1} \widehat{ \frac {1}{f_2} } (-s)& =
  (N+1) \sum_{s=1} ^N \widehat {f_1} (s) \widehat{ \frac {1}{f_2} } (-s) +\\
  &+  (N+1) \sum_{s=1} ^N \widehat {f_1} (s) \widehat{ \frac {1}{f_2} } (-s)
  \left( (1- \frac{s}{N} )^{\alpha_2+1}-1 \right)
  \end{align*}
   ce qui donne encore 
    \begin{align*} 
   (N+1) \sum _{s=1} ^N \widehat {f_1} (s) (1- \frac{s}{N}) ^
  {\alpha_2 +1} \widehat{ \frac {1}{f_2} } (-s)& =
  (N+1) \sum_{s=1} ^{+ \infty} \widehat {f_1} (s) \widehat{ \frac {1}{f_2} } (-s)
- (N+1) \sum_{s=N} ^{+\infty} \widehat {f_1} (s) \widehat{ \frac {1}{f_2} } (-s)\\
  &+  (N+1) \sum_{s=1} ^N \widehat {f_1} (s) \widehat{ \frac {1}{f_2} } 
  \left( (1- \frac{s}{N} )^{\alpha_2+1}-1 \right) 
  \end{align*}
  Nous pouvons donc conclure 
  \begin{eqnarray*}
 & & (N+1) \left(\hat{f_{1}} (0)\widehat{\left(\frac{1}{f_{2}}\right)} (0)+ 2 \sum _{s=1} ^N \widehat {f_1} (s) (1- \frac{s}{N}) ^
  {\alpha_2 +1} \widehat{ \frac {1}{f_2} } (-s)\right)\\
  &= &\trace  \left(T_N \left( \frac{f_1}{f_2} \right)\right) + 2 N^{2\alpha_2-2\alpha_1} 
  C_1 \left( \int_0 ^1 u^{2\alpha_2- 2 \alpha_1 -2} 
  \left( (1-u)^{\alpha_2+1} -1\right) du\right. -\\
  &-&\left.\int_1^{+ \infty} u^{2\alpha_2- 2 \alpha_1 -2} du\right)+o(N^{2\alpha_{2}-2\alpha_{1}})
  \end{eqnarray*}
     b) $ 0 > \alpha_2 - \alpha_1 >\frac{-1}{2}.$
On crit alors :
 \begin{align*} 
 & (N+1) \sum _{s=1} ^{N^{\prime}} \widehat {f_1} (s) (1- \frac{s}{N}) ^
  {\alpha_2 +1} \widehat{ f_{2}^{-1} } (-s) =
  (N+1) \sum_{s=1} ^{N^{\prime}} \widehat {f_1} (s) \widehat{ f_{2}^{-1} } (-s)
  -(\alpha_2+1)  \sum_{s=1} ^{N^{\prime}} s \widehat {f_1} (s) \widehat{ f_{2}^{-1}}(-s)  \\
  &+  (N+1) \sum_{s=1} ^{N^{\prime}} \widehat {f_1} (s) \widehat{f_{2}^{-1} } (-s)
  \left( (1- \frac{s}{N} )^{\alpha_2+1}-1+(\alpha_2+1)\frac{s}{N} \right) \\
  &=  (N +1)\sum_{s=1} ^{+\infty} \widehat {f_1} (s) \widehat{ f_{2}^{-1} } (-s)
  -(\alpha_2+1)  \sum_{s=1} ^{+\infty} s \widehat {f_1} (s) \widehat{ f_{2}^{-1} } (-s)-\\
  &- (N+1) \sum_{s=N+1} ^{+\infty} \widehat {f_1} (s) \widehat{ f_{2}^{-1} } (-s)
  +(\alpha_2+1)  \sum_{s=N+1} ^{+\infty} s \widehat {f_1} (s) \widehat{ f_{2}^{-1} } (-s) +\\
  &+ ( N +1) \sum_{s=1} ^{N^{\prime}} \widehat {f_1} (s) \widehat{ f_{2}^{-1} } 
  \left( (1- \frac{s}{N} )^{\alpha_2+1}-1+(\alpha_2+1)\frac{s}{N} \right) 
  \end{align*}
  D'o\`u, 
  \begin{eqnarray*}
  & & (N+1) \left(\widehat {f_{1}} (0) \widehat{\left(f_{2}^{-1}\right)} (0)+ 2 \sum _{s=1} ^{N^{\prime}} \widehat {f_1} (s) (1- \frac{s}{N}) ^
  {\alpha_2 +1} \widehat{ f_{2}^{-1}} (-s)\right)=\\
&=& \trace \left(T_N \left( \frac{f_1}{f_2} \right) \right) - (\alpha_2 +1)
 \langle f_{1} \vert f_{2}^{-1} \rangle _{2,1/2}
+\\
&+& 2 N^{2\alpha_2-2\alpha_1} C_1 \left( \int_ 0^1 u^{2\alpha_2-2 \alpha_1 -2}
 \left( (1-u)^{\alpha_2+1} -1 +(\alpha_2 +1) u \right) du \right.\\
& - & \left. \int_1 ^{+ \infty} u^{2\alpha_2 -2\alpha_1 -2} du +(\alpha_{2}+1)
  \int_1 ^{+ \infty} u^{2\alpha_2 -2\alpha_1 -1} du\right)+o(N^{2\alpha_{2}-2\alpha_{1}}).
  \end{eqnarray*}
 c)  $-\frac{1}{2} >\alpha_2 - \alpha_1 >-1$\\
 On utilise alors la dcomposition suivante 
 \begin{align*}
&  (N+1) \sum _{s=1} ^{N^{\prime}} \widehat {f_1} (s) \widehat{ f_{2}^{-1} } (-s)
  (1- \frac{s}{N}) ^ {\alpha_2 +1}  =\\
&  (N+1) \sum_{s=1} ^{N^{\prime}} \widehat {f_1} (s) \widehat{ f_{2}^{-1} } (-s)
  -(\alpha_2+1)  \sum_{s=1} ^{N^{\prime}} s \widehat {f_1} (s) \widehat{ f_{2}^{-1} } (-s)+
 \frac{(\alpha_2+1) \alpha_2}{2N} \sum_{s=1}^ {N^{\prime}}
s^2 \widehat {f_1} (s) \widehat{ \frac {1}{f_2} } (-s)+\\
& (N+1) \sum _{s=1} ^{N^{\prime}} \widehat {f_1} (s) \widehat{ f_{2}^{-1} } (-s)
  \left( (1- \frac{s}{N}) ^ {\alpha_2 +1} -1 
 + (\alpha_2+1) \frac{s}{N} 
  -\frac{\alpha_2(\alpha_2+1) s^2}{2 N^2} \right)
  \end{align*}
  Ce qui nous permet finalement d'\'ecrire,en utilisant les m\^{e}m`es calculs que ci-dessus  
  \begin{eqnarray*}
   & & (N+1) \left(\hat{f_{1}} (0) \widehat{\left(f_{2}^{-1}\right)} (0)+ 2 \sum _{s=1} ^{N^{\prime}} \widehat {f_1} (s) (1- \frac{s}{N}) ^
  {\alpha_2 +1} \widehat{ f_{2}^{-1}} (-s)\right)=\\
  &=& \trace\left( T_N \left( \frac{f_1}{f_2} \right)\right)  -  (\alpha_2+1) \langle f_{1} \vert f_{2}^{-1} \rangle _{2,1/2}
+ O(\frac{1}{N}).
\end{eqnarray*}
\subsubsection{Calcul de 
$ \displaystyle{
2 \Re\left(\sum_{s=1}^{N^{\prime}} \widehat{f_{1}} (s) (1- \frac{s}{N} ) ^{\alpha_{2}}
 \tau_{2} (\frac{s}{N}) \displaystyle{\left( \sum _{u\ge 0} u \beta_{u}^{(\alpha_{2})} 
\overline{\beta_{u+s}^{(\alpha_{2})}}\right)}\right)}
$
}
Pour all\'eger les calculs nous poserons dans la suite de cette d\'emonstration 
$$\left( \sum _{u\ge 0} u \beta_{u}^{(\alpha_{2})} 
\overline{\beta_{u+s}^{\alpha_{2}}}\right)
= \Sigma (s).$$
 Rappelons le r\'esultat, pour $s$ suffisamment grand  
 $$ \Sigma (s) =
 C_{2} s^{2 \alpha_{2}}, \quad \mathrm {avec} \quad C_{2}= \frac{C_{1}}{\alpha_{2}} .$$
 
 L\`a aussi nous allons devoir distinguer trois cas suivant les valeurs de $ \alpha_{2}- \alpha_{1}$.\\
 
  a) $\frac{1}{2}>\alpha_{2 }- \alpha_{1} >0$\\
 Dans ce cas la formule d'Euler et Mac-Laurin permet d'\'ecrire directement 
 \begin{eqnarray*}
& & 2 \Re\left(\sum_{s=1}^{N^{\prime}} \widehat{f_{1}} (s) (1- \frac{s}{N} ) ^{\alpha_{2}}
 \tau_{2} (\frac{s}{N}) \displaystyle{\left( \sum _{u\ge 0} u \beta_{u}^{(\alpha_{2})} 
\overline{\beta_{u+s}^{(\alpha_{2})}}\right)}\right)
=\\
 &=& 2 N^{2 \alpha_{2}- 2 \alpha_{1}}  \int _{0}^1 u^{2\alpha_{2}- 2\alpha_{1}-1 } (1-u) ^{\alpha_{2}}
 \tau_{2}(u) du +o(N^{2\alpha_{2}-2\alpha_{1}})
 \end{eqnarray*}

 b) $0> \alpha_{2}- \alpha_{1} > -\frac{1}{2}$\\
Nous allons utiliser la d\'ecomposition suivante 
 \begin{eqnarray*}
 & & \sum_{s=1}^{N^{\prime}} \widehat{f_{1}} (s) (1- \frac{s}{N} ) ^{\alpha_{2}}
 \alpha_{2}(2-\frac{s}{N}) \Sigma (s) =\\
&=& 2 \alpha_{2}  \sum_{s=1}^{N^{\prime}} \widehat{f_{1}} (s) (1- \frac{s}{N} ) ^{\alpha_{2}}
\Sigma (s) -
\frac {\alpha_{2}}{N} \sum_{s=1}^{N^{\prime}} s \widehat{f_{1}} (s) (1- \frac{s}{N} ) ^{\alpha_{2}}
\Sigma (s) \\
&=& 2 \alpha_{2} \left( \sum_{s=1}^{N^{\prime}} \widehat{f_{1}} (s) 
\Sigma (s) \left( (1- \frac{s}{N} )^{\alpha_{2}} -1 \right) +
  \sum_{s=1}^N \widehat{f_{1}} (s) 
\Sigma (s) \right)\\
&-& \frac{\alpha_{2}}{N}  \sum_{s=1}^{N^{\prime}} s \widehat{f_{1}} (s) 
\Sigma (s)
 (1-\frac{s}{N} )^{\alpha_{2}}. 
\end{eqnarray*}

    Ce qui donne, tous calculs faits 
 \begin{eqnarray*}
 & & \Re \left( \sum_{s=1}^{N^{\prime}} \widehat{f_{1}} (s) (1- \frac{s}{N} ) ^{\alpha_{2}}
\Sigma (s) \alpha_{2}( 2-\frac{s}{N})Ê\right) =
 2 \alpha_{2} \Re \left( \sum_{s=1}^{+ \infty} \widehat{f_{1}} (s) \Sigma (s) \right)\\
&+ &  C_{1}N^{2\alpha_{2}-2 \alpha _{1}} \left(2 \int_{0}^1u^{2 \alpha_{2}-2 \alpha_{1}-1} \left( (1-u) ^{\alpha_{2}} -1 \right)du\right.\\
&-& \left. \int_{0}^1u^{2 \alpha_{2}-2 \alpha_{1}}  (1-u) ^{\alpha_{2}} du 
- 2 \int_1^{+\infty}u^{2 \alpha_{2}-2 \alpha_{1} -1} du  \right) +o(N^{2\alpha_{2}-2\alpha_{1}}).
 \end{eqnarray*}
 On a de m\^{e}me 
  \begin{eqnarray*}
 &&  \Re\left(
\sum_{s=1}^{N^{\prime}} \widehat{f_{1}} (s) (1- \frac{s}{N} ) ^{\alpha_{2}+1}
2\Sigma (s)\right)=\\
&= &2\Re\left( \sum_{s=1}^{+\infty} \widehat{f_{1}} (s) 
\Sigma (s) \right)+\\
&+& 4 C_{1}N^{2\alpha_{2}-2 \alpha _{1}} \left( \int_{0}^1u^{2 \alpha_{2}-2 \alpha_{1}-1} \left( (1-u) ^{\alpha_{2}+1} -1 \right)du\right.\\
&-& 
\left. \int_1^{+\infty}u^{2 \alpha_{2}-2 \alpha_{1} -1} du  \right) +o(N^{2\alpha_{2}-2\alpha_{1}}).
 \end{eqnarray*}

Nous sommes donc ramen\'es \`a calculer 
$$4(\alpha_{2}+1) \Re(\displaystyle { \sum_{s=1} ^{+ \infty} \widehat {f_1} (s) \Sigma (s)}).$$
Pr\'ecisons cette quantit\'e quand les fonctions $f_{1}$ et $f_{2}$ sont paires.
Nous avons 
\begin{align*}
\sum_{u\ge 0} u \beta_{u} ^{(\alpha_{2})} \overline{\beta_{u+s} ^{(\alpha_{2})}}
&= - i\sum_{u\ge 0} u \widehat {g_{2}^{-1}} (u) \overline{\widehat {g_{2}^{-1}} (u+s)}\\
&=-i \langle \left( g_{2}^{-1}\right)' \vert \left( \chi^{-s} g_{2}^{-1}\right)\rangle 
= -i \langle \left( g_{2}^{-1}\right)' g_{2}\vert \left( \chi^{-s} f_{2}^{-1}\right)\rangle\\
&= \sum _{u \ge 0} u \widehat { \left( \ln g_{2}^{-1} \right) } (u) \overline{ \widehat{ \chi^{-s} f_{2}^{-1}} (u)}.
\end{align*}
Et de m\^{e}me nous obtenons :
\begin{align*}
&\sum_{u\ge 0} \overline{u \beta_{u} ^{(\alpha_{2})} } \beta_{u+s} ^{(\alpha_{2})}
=  \sum_{u\ge 0} \overline{ u \widehat {g_{2}^{-1}} (u)} \widehat {g_{2}^{-1}} (u+s)
\\
&= i \sum_{u\ge 0} \overline{\widehat {(g_{2}^{-1})'}} (u) \widehat {g_{2}^{-1} (u+s)}
= i \sum_{u\ge 0}\widehat { \overline{(g_{2}^{-1})'} }(-u) \overline         
 {\left(\overline{\widehat {g_{2}^{-1} }(u+s)}\right)}
\\
&= i \sum_{u\ge 0}\widehat { \overline{(g_{2}^{-1})'} }(-u)\overline{\widehat {\bar g_{2}^{-1}} (-u-s)}
=i \langle \left( \bar g_{2}^{-1}\right)' \vert \left( \chi^{s} \bar g_{2}^{-1}\right)\rangle 
\\
&= i \langle \left( \bar g_{2}^{-1}\right)' \bar g_{2}\vert \left( \chi^{s} f_{2}^{-1}\right)\rangle
= -\sum _{u \le 0} u \widehat { \left( \ln \bar g_{2}^{-1} \right) } (u) \overline{ \widehat{ \chi^{s} f_{2}^{-1}} (u)}.
\end{align*}

En remarquant que l'on a  $\pi_{+} \widehat {\ln (f_{2}^{-1})}(u) =  \widehat {\ln (g_{2}^{-1})}(u)$ si $u>0$
et $\pi_{-} \widehat {\ln (f_{2}^{-1})}(u) =  \widehat {\ln (\bar g_{2}^{-1})}(u)$ si $u< 0$ nous pouvons \'ecrire 
\begin{align*}
2 \Re \left( \sum _{s\ge 1} \widehat {f_{1}}(s) \Sigma (s)\right) &= 
\sum _{s\ge 1} \widehat {f_{1}}(s) \sum_{u \ge 0} u \widehat{ \ln f_{2}^{-1}} (u) \overline{ \widehat{ \chi^{-s} f^{-1}_{2}}(u)}
\\
& -\sum _{s\ge 1} \widehat {f_{1}}(-s) \sum_{u \le 0} u
\widehat{ \ln f_{2}^{-1}} (u) \overline{ \widehat{ \chi^{s} f^{-1}_{2}}(u)}.
\end{align*}

Nous avons ensuite 
\begin{align*}
\sum _{s\ge 1} \widehat{ (f_{1})}(s) \sum_{u \ge 0} u \widehat \ln f_{2}^{-1} (u) \overline{ \widehat{ \chi^{-s} f^{-1}_{2}}(u)}
&= \sum_{u \ge 0} u \widehat{ \ln f_{2}^{-1}} (u) \sum _{s\ge 1} \widehat{ (f_{1})}(s) \overline{ \widehat{ \chi^{-s} f^{-1}_{2}}(u)}\\
&= \sum_{u \ge 0} u \widehat {\ln f_{2}^{-1}} (u) \sum _{s\ge 1} \widehat {(f_{1})}(s) \overline{\widehat{f_{2}^{-1} }(u+s)}\\
&= \sum_{u \ge 0} u \widehat {\ln f_{2}^{-1}} (u) \sum _{s\ge 1} \widehat {(f_{1})}(s) \widehat{f_{2}^{-1} }(-u-s)
\end{align*}
et en remplaant $u$ par $-u$, en utilisant la parit de $f_2$
\begin{equation} \label{tortor}
\sum _{s\ge 1} \widehat{ (f_{1})}(s) \sum_{u \ge 0} u \widehat \ln f_{2}^{-1} (u) \overline{ \widehat{ \chi^{-s} f^{-1}_{2}}(u)}
=
\sum_{u \le 0} \vert u \vert \widehat {\ln f_{2}^{-1}} (u) \sum _{s\ge 1} \widehat {(f_{1})}(s) \widehat{f_{2}^{-1} }(u-s);
\end{equation}
en procdant de m\^{e}me il vient 
\begin{align*}
-\sum _{s\ge 1} \widehat (f_{1})(-s) \sum_{u \le 0} u \widehat{ \ln f_{2}^{-1} }(u) \overline{ \widehat{ \chi^{s} f^{-1}_{2}}(u)}
&= \sum_{u \le 0} \vert u \vert \widehat{ \ln f_{2}^{-1}} (u) \sum _{s\ge 1} \widehat {(f_{1})}(-s) \overline{ \widehat{ \chi^{-s} f^{-1}_{2}}(u)} \\
&= \sum_{u \le 0} \vert u \vert \widehat {\ln f_{2}^{-1}} (u) \sum _{s\ge 1} \widehat {(f_{1})}(-s) \overline{\widehat{f_{2}^{-1} }(u-s)}\\
&= \sum_{u \le 0} \vert u \vert \widehat {\ln f_{2}^{-1}} (u) \sum _{s\ge 1} \widehat {(f_{1})}(s) \widehat{f_{2}^{-1} }(s-u)
\end{align*}
ce qui donne, toujours gr\^{a}ce  la parit de $f_2$ :
\begin{equation} \label{tortor2}
-\sum _{s\ge 1} \widehat (f_{1})(-s) \sum_{u \le 0} u \widehat{ \ln f_{2}^{-1} }(u) \overline{ \widehat{ \chi^{s} f^{-1}_{2}}(u)}
=
 \sum_{u \le 0} \vert u \vert \widehat {\ln f_{2}^{-1}} (u) \sum _{s\ge 1} \widehat {(f_{1})}(s) \widehat{f_{2}^{-1} }(u-s).
\end{equation}
Nous pouvons donc crire 
\begin{equation} \label{tortor3}
2 \Re \left( \sum _{s\ge 1} \widehat {(f_{1})}(s) \Sigma (s)\right)
=
 2 \sum_{u \le 0} \vert u \vert \widehat {\ln f_{2}^{-1}} (u) \sum _{s\ge 1} \widehat {(f_{1})}(s) \widehat{f_{2}^{-1} }(u-s).
 \end{equation}
Considrons maintenant l'galit 
\begin {equation} \label{tortor4}
2 \left(\langle \ln f_2^{-1} \vert \frac{f_1}{f_2} \rangle_{2,1/2} 
- \langle \ln f_2^{-1} \vert \frac{1}{f_2} \rangle_{2,1/2} \widehat {f_{1}}(0)\right)
=+S_1+S_2+S_3+S_{4}
\end{equation}
avec 
$$S_1 =2\sum_{u \le 0} \vert u \vert \widehat {\ln f_{2}^{-1}} (u) \sum _{s\le -1} \widehat {(f_{1})}(s) \widehat{f_{2}^{-1} }(u-s)$$
$$S_2 =2
\sum_{u \ge 0} \vert u \vert \widehat {\ln f_{2}^{-1}} (u) \sum _{s\le -1} \widehat {(f_{1})}(s) \widehat{f_{2}^{-1} }(u-s)
 $$
 $$S_3 =2
 \sum_{u \ge 0} \vert u \vert \widehat {\ln f_{2}^{-1}} (u) 
 \sum _{s\ge 1} \widehat {(f_{1})}(s) \widehat{f_{2}^{-1} }(u-s),$$
 $$
 S_{4}= 2 
  \sum_{u \le 0} \vert u \vert \widehat {\ln f_{2}^{-1}} (u) 
 \sum _{s\ge 1} \widehat {(f_{1})}(s) \widehat{f_{2}^{-1} }(u-s).$$
 
 Nous avons clairement 
 $ S_4=S_2 = \Re \left( \displaystyle{\sum _{s\ge 1} \widehat {(f_{1})}(s) \Sigma (s)}\right)$
 et 
$S_1 =S_3$.
Nous avons d'autre part 
\begin{align*}
S_1 &= \sum_{u \le 0} \vert u \vert \widehat {\ln f_{2}^{-1}} (u) \sum _{s\le -1} \widehat {(f_{1})}(s) \widehat{f_{2}^{-1} }(u-s)\\
&= \sum_{u \le 0} \vert u \vert \widehat {\ln f_{2}^{-1}} (u) \sum _{s\le -1} \widehat {(f_{1})}(s) \overline{\widehat{f_{2}^{-1} }(s-u)}\\
&= \sum_{u \le 0} \vert u \vert \widehat {\ln f_{2}^{-1}} (u) \sum _{s\ge 1} \widehat {(f_{1})}(s)\overline{ \widehat{f_{2}^{-1} }(s+u)}\\
&= -\sum _{s\ge 1} \widehat {(f_{1})}(s)\sum_{u \le 0}  u \widehat {\ln f_{2}^{-1}} (u) \overline{ \widehat{ \chi ^{-s}f_{2}^{-1} }(u)}.
\end{align*}
Nous pouvons maintenant constater que 
\begin{eqnarray*}
-\sum_{u \le 0}  u \widehat {\ln f_{2}^{-1}} (u) \overline{ \widehat{ \chi ^{-s}f_{2}^{-1} }(u)} &= &
-\sum_{u \le 0} \widehat u {\ln \bar g_{2}^{-1}} (u) \overline{ \widehat{ \chi ^{-s}f_{2}^{-1} }(u)}
\\
&=& i \sum_{u \le 0} \widehat { \left( \ln \bar g_{2}^{-1}\right)'} (u)
\overline{ \widehat{ \chi ^{-s}f_{2}^{-1} }(u)}\\
&=& i \langle \left( \ln \bar g_{2}^{-1}\right)' \vert \chi ^{-s}f_{2}^{-1}
\rangle = i \langle \frac{\left( \bar g_{2}^{-1}\right)'}{\bar g_2^{-1}} \vert \chi ^{-s}f_{2}^{-1}\rangle \\
&=& i \langle \left( \bar g_{2}^{-1}\right)' \vert \chi ^{-s} \bar g _{2}^{-1} \rangle = i \sum_{u\le 0} \widehat {\left( \bar g_{2}^{-1}\right)'}
(u) \overline{\widehat{\chi ^{-s} \bar g _{2}^{-1} }(u)}\\
&=&
- \sum_{u\le 0} u \widehat{\bar g_2^{-1}}(u)  
\widehat{\chi ^{s}  g _{2}^{-1} }(-u)
= - \sum_{u \le -s } u \overline{\beta_{-u}^{(\alpha_2)} } 
 \beta_{-u-s}^{(\alpha_2)}
\end{eqnarray*}
 Il est maintenant facile de se rendre compte que 
 $$
 - \sum_{u \le -s } u \overline{\beta_{-u}^{(\alpha_2)} } 
 \beta_{-u-s}^{(\alpha_2)} =
 \sum_{v \ge s } v \overline{\beta_{v}^{(\alpha_2)} } 
 \beta_{v-s}^{(\alpha_2)}
 = \sum_{w\ge 0 } (w+s) \overline{\beta_{w+s}^{(\alpha_2)} } 
 \beta_{w}^{(\alpha_2)}
 $$
 et puisque 
 $$
 s \sum_{w\le 0 }  \overline{\beta_{w+s}^{(\alpha_2)} } 
 \beta_{w}^{(\alpha_2)} = 
 s \langle g^{-1}_2 \vert \chi ^s g^{-1}_2 \rangle 
= s \langle 1 \vert \chi ^s f^{-1}_2 \rangle 
 = s\widehat {f_2} (-s)
 $$
 nous pouvons finalement conclure 
 \begin {equation} \label{tortor 6}
2 \left(\langle \ln f_2^{-1} \vert \frac{f_1}{f_2} \rangle_{2,1/2} 
- \langle \ln f_2^{-1} \vert \frac{1}{f_2} \rangle_{2,1/2} \widehat {f_{1}} (0) \right)
= 4 \Re \left( \sum _{s\ge 1} \widehat {(f_{1})}(s) \Sigma (s)\right)+ 2 \langle f_1 \vert f_2^{-1}\rangle_{2,1/2}.
\end{equation}
 Autrement dit 
 \begin {equation} \label{tortor 7}
 \Re \left( \sum _{s\ge 1} \widehat {(f_{1})}(s) \Sigma (s)\right) 
 = \frac{1}{2} \left( \langle \ln f_2^{-1} \vert \frac{f_1}{f_2} \rangle_{2,1/2} \widehat {f_{1} }(0) 
- \langle \ln f_2^{-1} \vert f_2^{-1} \rangle_{2,1/2}-\langle f_1 \vert f_2^{-1}\rangle_{2,1/2}\right)
\end{equation}
 Nous pouvons donc conclure que 
  $$ 
2 \Re\left(\sum_{s=1}^N \widehat{f_{1}} (s) (1- \frac{s}{N} ) ^{\alpha_{2}}
 \tau_{2} (\frac{s}{N}) \left( \sum _{u\ge 0} u \beta_{u}^{(\alpha_{2})} 
\overline{\beta_{u+s}^{(\alpha_{2})}}\right)\right)
$$
vaut 
$$ -2(\alpha_2+1)  \left( \langle \ln f_2^{-1} \vert \frac{f_1}{f_2} \rangle_{2,1/2} 
- \langle \ln f_2^{-1} \vert f_2^{-1} \rangle_{2,1/2} \widehat {f_{1} }(0) -\langle f_1 \vert f_2^{-1}\rangle_{2,1/2}\right).$$
  c) $ -\frac{1}{2} > \alpha_{2}-\alpha_{1}> -1$ 
 Nous obtenons alors 
  \begin{eqnarray*}
 & & \sum_{s=1}^{N^{\prime}} \widehat{f_{1}} (s) (1- \frac{s}{N} ) ^{\alpha_{2}}
\left( \sum _{u\ge 0} u \beta_{u}^{(\alpha_{2})} 
\overline{\beta_{u+s}^{(\alpha_{2})}}\right) \alpha_{2}(2-\frac{s}{N})=\\
&=& 2 \alpha_{2} \left( \sum_{s=1} ^{N^{\prime}} \widehat {f_{1}}(s) 
\left( \sum _{u\ge 0} u \beta_{u}^{(\alpha_{2})} 
\overline{\beta_{u+s}^{(\alpha_{2})}}\right)
 \left( (1- \frac{s}{N} )^{\alpha_{2}} -1 + \alpha_{2} \frac{s}{N}\right) 
 +  \sum_{s=1} ^N \widehat {f_{1}}(s) 
\left( \sum _{u\ge 0} u \beta_{u}^{(\alpha_{2})} 
\overline{\beta_{u+s}^{(\alpha_{2})}}\right)\right. -\\
 &-&\left. \frac{\alpha_{2}}{N} \sum_{s=1}^{N^{\prime}} s \widehat {f_{1}} (s)
 \left( \sum _{u\ge 0} u \beta_{u}^{(\alpha_{2})} 
\overline{\beta_{u+s}^{(\alpha_{2})}}\right) \right)\\
 &-& \frac{\alpha_{2}} {N}  \left( \sum_{s=1} ^{N^{\prime}}  s \widehat {f_{1}}(s) 
 \left( \sum _{u\ge 0} u \beta_{u}^{(\alpha_{2})} 
\overline{\beta_{u+s}^{(\alpha_{2})}}\right)+ \sum_{s=1} ^N \widehat {f_{1}}(s) 
 \left( \sum _{u\ge 0} u \beta_{u}^{(\alpha_{2})} 
\overline{\beta_{u+s}^{(\alpha_{2})}}\right) \left( (1-\frac{s}{N} )^{\alpha_{2}}-1 \right)
 \right)
 \end{eqnarray*}
 et les r\'esultats du points pr\'ec\'edent s'appliquent imm\'ediatement.

   \subsubsection{ Calcul de $\displaystyle{N^{2\alpha_{2}} \sum_{s=0}^{N^{\prime}} (1-\frac{s}{N}) ^{2 \alpha_{2}} F_{1}(\frac{s}{N}) \widehat {f_{1}} (s) N^{2\alpha_{2}}}$}
  Ici seul le signe de $\alpha_{1}$ est d\'eterminant. (rappelons que $\alpha_{2}$ est n\'egatif)\\
  a)  $\alpha_{1}<0$\\
  Nous avons alors ais\'ement 
  \begin{eqnarray*}
 & & N^{2\alpha_{2}} \sum_{s=0}^{N^{\prime}} (1-\frac{s}{N} )^{2\alpha_{2}} F_{1}(\frac{s}{N}) \widehat F_{1} (s)  =\\
&=& N^{2\alpha_{2}-2 \alpha_{1}} \frac{c_{1} (1)} {\Gamma^2 (\alpha_{1} )}
\int _{0}^1 (1-t) ^{\alpha_{2} } F_{1}(t) t ^{-2\alpha_{1}-1} dt 
\end{eqnarray*}
 b)  $\alpha_{1}>0$\\
  Nous utilisons alors la d\'ecomposition
   \begin{eqnarray*}
   & & N^{2\alpha_{2} }\sum_{s=0}^{N^{\prime}} (1-\frac{s}{N}) ^{2\alpha_{2}} F_{1}(\frac{s}{N}) \widehat f_{1} (s)  =\\
 &=&  N^{2\alpha_{2} }\sum_{s=0}^{N^{\prime}} \left( F_{1}(\frac{s}{N}) (1-\frac{s}{N} )^{2\alpha_{2}} -F_{1}(0)\right)  \widehat f_{1} (s) 
 +  N^{2\alpha_{2} } F_{1}(0)  \sum_{s=0}^{N^{\prime}}  \widehat f_{1} (s) 
\end{eqnarray*}
qui donne l'\'egalit\'e suivante, en se souvenant que $\displaystyle{\sum_{-\infty}^{+ \infty} \widehat{f_{1}}(s) =0}$
\begin{eqnarray*}
   & & N^{2\alpha_{2} }\sum_{s=0}^N (1-\frac{s}{N} ^{\alpha_{2}} F_{1}(\frac{s}{N}) \widehat f_{1} (s)  =\\
&=& N^{2\alpha_{2}-2\alpha_{1}} \frac{c_{1} (1)} {\Gamma^2 (\alpha_{1} )} \left(
\int_{0}^1 \left( F_{1} (t) (1-t) ^{2\alpha_{2}} - F_{1}(0) \right) t^{-2\alpha_{1}-1} dt 
- F_{1}(0) K \right)
\end{eqnarray*}
avec 
$$ K = \int_{1}^{+\infty} t^{-2\alpha_{1}-1}dt  +
 \frac{\widehat f_(0)}{2}.$$
    \subsection {D\'emonstation du th\'eor\`eme dans le cas o\`u $\alpha_{2} \in ]0,\frac{1}{2}[$}
  \subsubsection {Calcul de 
  $ \displaystyle{N \sum _{s=1} ^N \widehat {f_1} (s) (1- \frac{s}{N}) ^
  {\alpha_2 +1} \widehat{ \frac {1}{f_2} } (-s) }$}
  Rappelons que dans ce cas 
  $$ \trace \left(T_{N}^{-1} (f_{2})\right)) = (N+1) \widehat{\left(\frac{1}{f_{2}}\right)} (0) +N^{2\alpha_{2} }
  \frac{2 K_{\alpha_{2}} }{c_{2} (1) \Gamma ^2 (\alpha_{2})}+ o(N^{2\alpha_{2}})$$
  avec 
  $$ K_{\alpha_{2}} = \int_{1/2}^1 \int_{0}^x t^{2\alpha_{2}-2}\left( (1-t)^{2\alpha_{2}}-1 \right) - 
  t^{2\alpha_{2}}(1-t)^{2\alpha_{2}-2} dt dx + \frac{1}{2\alpha_{2} (2\alpha_{2}-1)}.$$
  Comme dans le cas n\'egatif trois cas sont \`a distinguer.\\
a)$ \frac{1}{2} <\alpha_{2} - \alpha_{1} <1$.\\
Nous avons alors imm\'ediatement l'\'egalit\'e :
$$
 N \sum_{s=1}^N (1- \frac{s}{N}) ^{\alpha_{2}+1} \widehat{f_{1}}(s) \widehat {\frac{1}{f_{2}}} (-s) 
=
N^{2\alpha_{2}-2\alpha_{1}} C_{1}  \int _{0}^1(1-u)^{\alpha_{2}+1} u^{2\alpha_{2}-2\alpha_{1}-2} du 
+o(N^{2\alpha_{2}-2\alpha_{1}}).
$$
  b) $ 0 < \alpha_{2}-\alpha_{1}<\frac{1}{2}.$\\
  Avec les m\^{e}mes id\'ees que dans le cas n\'egatif il vient 
  \begin{eqnarray*}
  & & N \sum_{s=1}^N (1- \frac{s}{N}) ^{\alpha_{2}+1} \widehat{f_{1}}(s) \widehat {\frac{1}{f_{2}}} (-s) 
  =  N \sum_{s=1}^N  \widehat{f_{1}}(s) \widehat {\frac{1}{f_{2}} }(-s) 
  \left((1- \frac{s}{N}) ^{\alpha_{2}+1}-1\right) +N \sum_{s=1}^N  \widehat{f_{1}}(s) \widehat {\frac{1}{f_{2}} }(-s) =\\
  &=& \trace\left( T_{N}\left( \frac{f_{1}}{f_{2}}\right) \right)+ N^{2\alpha_{2}-2\alpha_{1}} C_{1}
  \left( \int_{0}^1 \left( (1-u)^{\alpha_{2}+1} -1\right) u^{2\alpha_{2}-2\alpha_{1}-2} du
-\int_{1}^{+ \infty} u^{2\alpha_{2}-2\alpha_{1}-2}du\right)\\
&+& o(N^{2\alpha_{2}-2\alpha_{1}}).
\end{eqnarray*}
c)  $- \frac{1}{2} <\alpha_{2} - \alpha_{1} < 0$.\\
Nous obtenons avec les m\^{e}mes justifications que dans 
le cas o $\alpha_2$ est ngatif  
 \begin{eqnarray*}
 & & N \sum_{s=1}^{N^{\prime}}  (1- \frac{s}{N}) ^{\alpha_{2}+1} \widehat{f_{1}}(s) \widehat {\frac{1}{f_{2}}} (-s) 
=  N \sum_{s=1}^{N^{\prime}}  \widehat{f_{1}}(s) \widehat {\frac{1}{f_{2}} }(-s) 
  \left((1- \frac{s}{N}) ^{\alpha_{2}+1}-1+ (\alpha_{2}+1) \frac{s}{N}\right)+\\
  &+& N \sum_{s=1}^{N^{\prime}}  \widehat{f_{1}}(s) \widehat {\frac{1}{f_{2}} }(-s) 
  - (\alpha_{2}+1)  \sum_{s=1}^{N^{\prime}}  s \widehat{f_{1}}(s) \widehat {\frac{1}{f_{2}} }(-s)
  \end{eqnarray*}
  D'o\`u $2 \displaystyle{\Re \left( N \sum _{s=1} ^{N^{\prime}} \widehat {f_1} (s) 
  (1- \frac{s}{N}) ^ {\alpha_2 +1} \widehat{ \frac {1}{f_2} } (-s) \right) + \widehat f_{1} (0)  \trace \left(T_{N}^{-1} (f_{2})\right))}$
  est \'egal \`a 
  \begin{eqnarray*}
  \trace\left( T_{N}\left( \frac{f_{1}}{f_{2}}\right) \right) 
 & -&(\alpha_{2}+1) \langle f_{1}\vert f_{2}^{-1}\rangle_{2,1/2} 
 + N^{2\alpha_{2} }
  \frac{2 K_{\alpha_{2}} }{c_{2} (1) \Gamma ^2 (\alpha_{2})}+\\
 & +& N^{2\alpha_{2}-2\alpha_{1}} C_{1} \left( \int _{0}^1 \left( (1-u)^{}{\alpha_{2}+1}-1+u (\alpha_{2}+1) \right) 
 u^{2\alpha_{1}-2\alpha_{1}-2} du \right.-\\
  &-& \left.\int_{1}^{+ \infty} u^{2\alpha_{2}-2\alpha_{1}-2} du + (\alpha_{2}+1) 
 \int_{1}^{+ \infty} u^{2\alpha_{1}-2\alpha_{1}-1}du \right).
 \end{eqnarray*}  
  \subsubsection{ Calcul de $\displaystyle{N^{2\alpha_{2}} \sum_{s=0}^{N^{\prime}} (1-\frac{s}{N}) 2^{\alpha_{2}} F_{2}(\frac{s}{N}) \widehat {f_{1}} (s) }$}
  Ce calcul se traite comme dans le cas pr\'ec\'edent.
  \section{D\'emonstration du th\'eor\`eme (\ref {thpuissance})}
  La d\'emonstration n\'ecessite deux lemmes. 
  \begin{lemme} \label{derivpuissance}
 Avec les hypoth\`eses du th\'eor\`eme (\ref{thpuissance}) f il existe une constante strictement positive t $K_{0}$ ind\'ependante de $s$ telle que pour 
 tout entier $s$ on ait 
 $$\Bigl\vert \trace \left( T_{N}(f_{1}) T_{N}^{-1} (f_{2})\right)^s - T_{N}\left( \left(\frac{f_{1}}{f_{2}} \right) ^s\right) \Bigr\vert \le 
 K_{0} s \left( 2 \Vert f_{1} \Vert _{\infty} \Vert f_{2}^{-1} \Vert _{\infty} \right)$$
  \end{lemme}
  \begin{lemme} \label{reste}
 Avec les hypoth\`eses du th\'eor\`eme (\ref{thpuissance})  n\'egatif posons 
 pour tout entier $N$ et $t \in \Delta$ 
 $$
 R_{1,N}(t) =  \trace \left( T_{N}(f_{1}) T_{N}^{-1} (f_{2})\right)- T_{N} \left(\frac{f_{1}}{f_{2}} \right) -\Psi_{1}(t).$$
 Alors pour tout $N$ la fonction $R_{1,N}$ est analytique sur un voisinage de z\'ero et pour tout entier naturel non nul $p$ 
 on a 
 $\displaystyle{\lim_{N \rightarrow + \infty }R_{1,N}^{(p)}(0)=0 }$.
  \end{lemme}
   \subsection{D\'emonstration du lemme (\ref {derivpuissance})}
 \begin{remarque}
 Dans la suite de la d\'emonstration nous noterons $\Vert  \Vert_{1}$ la norme 
 $\Vert A \Vert _{1} = (\trace (^t A A)^{1/2})$.
 Si $\Vert \Vert$ d\'esigne la norme classique des matrices, rappelons les propri\'et\'es bien connues 
 $$ \trace A \vert \le \Vert A \Vert _{1},\quad
  \Vert A B\Vert_{1} \le \Vert A \Vert_{1} \Vert B \Vert.$$
 \end{remarque}
 \begin{remarque}
 On rappelle \'egalement que si $\rho$ est un endomorphisme sym\'etrique r\'eelle dans un espace de dimension n
 on a, si $ \lambda_{1}, \lambda_{2},\cdots, \lambda_{n} $ sont les valeurs propres de $\rho$ :
 $ \Vert \rho \Vert = \max \{ \lambda_{i} / 1 \le i \le n \}$.
 \end {remarque} 
 \begin{remarque}
 Enfin on utilisera que si $f$ est une fonction positive sur $]-Ê\pi, \pi]$ qui ne s'annule qu'en un nombre fini de 
 points l'op\'erateur $T_{N}(f^{-1}) - T_{N}^{-1}(f)$ est un op\'erateur positif.
 \end{remarque}
 Nous devons d'abord rappeler les r\'esultats techniques suivants 
 \begin{lemme} \label {norme1}
 Si $h_{1}$ et $h_{2}$ sont dans $A(\mathbb T, 1/2)$ on a 
 $$ \Vert T_{N}(h_{1} h_{2}) -T_{N}(h_{1}) T_{N}(h_{2}) \Vert _{1} \le \Vert h_{1}\Vert _{2,1/2} \Vert h_{2}\Vert _{2,1/2}.$$
  \end{lemme}
  et \'egalement 
  \begin{lemme} \label {norme12}
Il existe une constante $C>0$ telle que si $h_{1}$ et $h_{2}$ sont dans $A(\mathbb T, 1/2)$ on ait  
  $$ \Vert h_{1} h_{2} \Vert_{2,1/2} \le C (\Vert h_{1} \Vert_{\infty} +\Vert h_{1} \Vert _{2,1/2}) 
   (\Vert h_{2} \Vert_{\infty} +\Vert h_{2} \Vert _{2,1/2}).$$
  \end{lemme}
  Ces deux r\'esultats peuvent se trouver dans \cite {Seg}.\\
 Ecrivons 
 \begin{align*}
  \trace \left( T_{N} (f_{1}) T_{N}^{-1} (f_{2})\right)^s - T_{N}\left( \left(\frac{f_{1}}{f_{2}} \right) ^s\right) &=&
  \trace \left(\left( T_{N} (f_{1}) T_{N}^{-1} (f_{2})\right)^s - \left(T_{N} \left(\frac{f_{1}}{f_{2}} \right)\right)^s\right) +\\
  + \trace\left( \left(T_{N} \left(\frac{f_{1}}{f_{2}} \right)\right)^s- T_{N}\left( \left(\frac{f_{1}}{f_{2}} \right) ^s\right)\right).
  \end{align*}
  En utilisant 
  $$ \trace (A^s -B^s) = \trace \left( (A-B) (A^{s-1} + A^{s-2} B + \cdots +B^{s-1} \right)$$
  nous pouvons \'ecrire 
$$
\Bigl \vert \trace \left(\left( T_{N} (f_{1}) T_{N}^{-1} (f_{2})\right)^s - \left(T_{N} \frac{f_{1}}{f_{2}} \right)^s\right)\Bigr \vert  \le
 S_{1} S_{2}$$
 En posant 
 $S_{1}= \Vert T_{N} (f_{1}) T_{N}^{-1} (f_{2})- T_{N} \frac{f_{1}}{f_{2}}\Vert_{1}$
 et  
 $ S_{2}= \Vert \sum_{k=1}^{s-1} \left (T_{N} (f_{1}) T_{N}^{-1} (f_{2})\right)^k \left( T_{N}\frac {f_{1}}{f_{2}}\right)^{s-1-k} \Vert. $\\
Il vient alors, en utilisant \cite{RS09}, plus le fait que 
 $ T_{N} (f_{2}^{-1})-T_{N}^{-1}(f_{2}) $ est un op\'erateur positif et aussi le lemme \ref{norme1}
 \begin{align*}
 S_{1} &\le \Vert T_{N}(f_{1}) \Vert \Vert T_{N}^{-1}(f_{2})  - T_{N} (f_{2}^{-1}) \Vert_{1}
+ \Vert T_{N}(f_{1}) T_{N}(f_{2}^{-1}) - T_{N} \frac{f_{1}} {f_{2}}\Vert_{1}\\
&\le  O(N^{2\alpha_{2}}) + \Vert f_{1}\Vert _{2,1/2} \Vert f_{2}^{-1}\Vert _{2,1/2}.
\end{align*}
D'autre part nous obtenons facilement 
\begin{align*}
S_{2} &\le \sum_{p=1} ^s \Vert f_{1} \Vert _{\infty}^p \Vert \left(T_{N}(f_{2})\right)^{-1} \Vert ^p 
\Vert \frac{f_{1}}{f_{2}} \Vert _{\infty}^{s-1-p}
 \le \sum_{p=1} ^s \left( \Vert f_{1} \Vert _{\infty} \Vert f_{2}^{-1} \Vert_{\infty} \right)^p \Vert \frac{f_{1}}{f_{2}} \Vert _{\infty}^{s-1-p}\\
&\le  \sum_{p=1} ^s  \left( \Vert f_{1} \Vert _{\infty} \Vert f_{2}^{-1} \Vert_{\infty} \right)^{s-1}
 \le (s-1)   \left( \Vert f_{1} \Vert _{\infty} \Vert f_{2}^{-1} \Vert_{\infty} \right)^{s-1}
\end{align*}
Etudions maintenant 
$
\trace\left( \left(T_{N} \left(\frac{f_{1}}{f_{2}} \right)\right)^s- T_{N}\left( \left(\frac{f_{1}}{f_{2}} \right) ^s\right)\right)
$.
On a 
$$
\trace\left( \left(T_{N} \left(\frac{f_{1}}{f_{2}} \right)\right)^s- T_{N}\left( \left(\frac{f_{1}}{f_{2}} \right) ^s\right)\right)
\le \Bigl\Vert  \left(T_{N} \left(\frac{f_{1}}{f_{2}} \right)\right)^s- T_{N}\left( \left(\frac{f_{1}}{f_{2}} \right) ^s\right) \Bigr \Vert _{1}
$$
et on pose 
$$
G_{s}=\Bigl \Vert  \left(T_{N} \left(\frac{f_{1}}{f_{2}} \right)\right)^s- T_{N}\left( \left(\frac{f_{1}}{f_{2}} \right) ^s\right) \Bigr \Vert _{1}
 $$
 On a 
 \begin{eqnarray*} 
G_{s} &\le&\Bigl \Vert  \left( T_{N} \frac{f_{1}}{f_{2}} T_{N} \frac{f_{1}}{f_{2}} \right)^{s-1} 
 -   T_{N} \frac{f_{1}}{f_{2}} \left( T_{N} \left(\frac{f_{1}}{f_{2}} \right) ^{s-1}\right) \Bigr\Vert _{1} 
 \\
& +& \Bigl  \Vert T_{N} \frac{f_{1}}{f_{2}} \left(T_{N}\left( \frac{f_{1}}{f_{2}} \right) ^{s-1}\right)
-\left(  T_{N}\left(\frac {f_{1}}{f_{2}} \right) ^s\right)\Bigr  \Vert _{1}.
\end{eqnarray*}
Et en utilisant le lemme (\ref{norme1}) nous obtenons 
$$\Bigl\Vert T_{N} \frac{f_{1}}{f_{2}} \left(T_{N} \left(\frac{f_{1}}{f_{2}} \right) ^{s-1}\right)
- \left( T_{N}\left(\frac{f_{1}}{f_{2}} \right) ^s\right) \Bigr\Vert _{1}
\le 
 \Bigl \Vert \frac{f_{1}}{f_{2}} \Bigr\Vert _{2,1/2} \Bigl\Vert \left( \frac{f_{1}}{f_{2}} \right)^{s-1}\Bigr\Vert _{2,1/2}
 $$
 ce qui d'apr\`es le lemme \ref{norme12} s'\'ecrit aussi 
  $$\Bigl\Vert T_{N} \frac{f_{1}}{f_{2}} \left(T_{N} \left(\frac{f_{1}}{f_{2}} \right) ^{s-1}\right)
- \left( T_{N}\left(\frac{f_{1}}{f_{2}} \right) ^s\right) \Bigr\Vert _{1}
\le
  C \Bigl \Vert \frac{f_{1}}{f_{2}} \Bigr\Vert _{2,1/2} \left( \Bigl \Vert \frac{f_{1}}{f_{2}} \Bigr\Vert _{2,1/2}
 + \Bigl \Vert \frac{f_{1}}{f_{2}} \Bigr\Vert _{\infty} \right)^{s-1}.
 $$
 En posant $M =  \Bigl \Vert \frac{f_{1}}{f_{2}} \Bigr\Vert _{2,1/2}
 + \Bigl \Vert \frac{f_{1}}{f_{2}} \Bigr\Vert _{\infty}$ nous avons l'in\'egalit\'e 
 $$G_{s}\le  M G_{s-1}+  C M^s. $$
Ce qui donne finalement 
$$ G_{s} \le C   s  M^s.$$
D'o\`u le lemme( \ref {derivpuissance}).
 \subsection{D\'emonstration du lemme (\ref{reste})}
Un d\'eveloppement en s\'erie enti\`ere permet d'\'ecrire, pour tout $t \in \Delta$, 
$$
R_{1,N}(t) = \trace \left( T_{N} (f_{1}) \left( T_{N}  f_{2,t} \right)^{-1} - T_{N}\left( \frac{f_{1}}{f_{2,t}}\right)\right)
-\Psi_{1}(t) $$
soit 
$$R_{1,N}(t) = \sum_{s=0} ^{\infty} t^s  (-1)^s  \trace \left(  T_{N} (f_{1}) T_{N}^{-1} (f_{2})\right) ^{s+1} -\trace  T_{N}
 \left( \left(\frac{f_{1}}{f_{2}}\right)^{s+1}\right)  - \sum_{s=0} ^{\infty} t^s \frac{\Psi_{1}^{(s)} (0) }{s!} 
 $$
Ce qui donne, pour tout entier naturel $q$ (en d\'erivant) 
\begin{align*}
\Bigl \vert R_{1,N}^{(q)} (t) - R_{1,N}^{(q)} (0) \Bigr \vert &= \sum_{s\ge q+1} (-1)^s s \cdots (s-q+1) t^{s-q} 
\trace \left(  T_{N} (f_{1}) T_{N}^{-1} (f_{2})\right) ^{s+1} - 
\trace T_{N}
 \left( \left(\frac{f_{1}}{f_{2}}\right)^{s+1}\right)\\ 
 &- \sum_{s\ge q+1} t^{s-q}  \frac{ \Psi_{1}^{(s)} (0) }{(s-p)!}. 
 \end{align*}
nous obtenons alors, en utilisant le lemme (\ref{derivpuissance}
$$ \Bigl \vert R_{1,N}^{(q)} (t) - R_{1,N}^{(q)} (0) \Bigr \vert \le \vert t \vert \left( \sum_{s \ge p+1} \vert t^{s-p+1}\vert  K_{0}^s s \cdots (s-q+1)
+ \sum_{s \ge p+1} \vert t^{s-p+1} \vert  \frac{\vert \Psi _{1}^{(s)} (0)\vert } {(s-q)!} \right) $$
et 
$$\Bigl \vert R_{1,N}^{(q)} (t) - R_{1,N}^{(q)} (0) \Bigr \vert \le \vert t \vert K_{0}^{p+1}  \phi_q (t)$$
o la fonction $\phi_q$ est borne sur un voisinage de zro. Ce r\'esultat peut encore s'\'ecrire, si $p$ est comme dans l'\'enonc\'e du lemme 
\begin{equation} 
\label{formule} \forall \epsilon>0, \quad \forall q \in \mathbb N^* \quad \exists \delta _{\epsilon,q} \quad \mathrm{t. q.} \quad 
\forall N \quad \Bigl \vert R_{1,N}^{(q)} (t) - R_{1,N}^{(q)} (0) \Bigr \vert  \le \epsilon.
\end{equation}
Soit maintenant $t$ fix dans $]-\delta _{\epsilon,1},  \delta _{\epsilon,1}[$. Pour tout entier $N$ il existe $t_{N}$ compris entre 
$0$ et $t$ tel que 
$ R_{1,N}(t) = R_{1,N}(0) +  t R'_{1,N}(t_{N})$. C'est une cons\'equence du th\'eor\`eme \ref{thprincipal} que 
la limite des suites $ R_{N}(0)$ et $R_{1,N}(t)$ est nulle. On a donc $ \displaystyle{ \lim_{N \rightarrow + \infty } t R'_{1,N}(t_{N})} =0$ et 
$ \displaystyle{ \lim_{N \rightarrow + \infty }  R'_{1,N}(t_{N})} =0$.
Puisque 
$t_{N} \in ] - \delta _{\epsilon,1}, \delta _{\epsilon,1}[$  l'in\'egalit\'e 
(\ref{formule}) permet de conclure $ \displaystyle{ \lim_{N \rightarrow + \infty }R'_{1,N}(t) =0}$.\\
Supposons maintenant la propri\'et\'e obtenue pour tout entier $q$ strictement inf\'erieur \`a $p$. 
Soit cette fois $t$ dans $]-\delta _{\epsilon,p}, \delta _{\epsilon,p}[$. Pour tout entier $N$ la formule de Taylor-Lagrange nous donne 
l'existence d'au moins un r\'eel $t_{N}$ tel que 
$$R_{1,N} \sum_{j=0}^{p-1} \frac {t^j}{j!} R_{1,N}^{(s)} (0) + \frac{t^p} {p!} R_{1,N}^{(p)} (t_{N}).$$
Le th\'eor\`eme \ref{thprincipal} et l'hypoth\`ese de r\'ecurrence donnent de la m\^{e}me mani\`ere que ci-dessus 
$ \displaystyle{ \lim_{N \rightarrow + \infty }  R^{(p)}_{1,N}(t_{N})} =0$ et, avec  
 $t_{N} \in ] - \delta _{\epsilon,p}, \delta _{\epsilon,p}[$ l'ingalit (\ref{formule}) donne 
  $ \displaystyle{ \lim_{N \rightarrow + \infty }  R^{(p)}_{1,N}(0)} =0$.
  \subsection{Fin de la d\'emonstration du th\'eor\`eme \ref{thpuissance}}
  Le th\'eor\`eme s'obtient alors facilement pour $\alpha_{2}<0$ \`a partir du lemme \ref{reste} et en d\'erivant suffisamment de fois en zro la formule 
  $$ \Tr  \left( T_{N} (f_{1}) T_{N}^{-1} (f_{2,t}) \right) - \trace \left(T_{N} \frac {f_{1}}{f_{2,t}}\right) = \Psi_{1}(t) +R_{1,N}(t).$$
    \subsection{Dmonstration du corollaire \ref{corollaire}}
  Posons $\tilde h_2 =h_2^{-1}$. 
  Posons 
  $$ T_N h_1 T_N^{-1} \tilde h_2=A \quad  T_N (h_1) T_N \frac{1}{\tilde h_2}=B.$$
 Pour tout entier naturel $s$ nous avons
  $$ \trace (A^s -B^s ) = \trace\left( (A-B) \left(  \sum _{j=0}^{s-1}
  A^{j} B^{s-1-j}\right)\right)$$
  ce qui permet d'crire 
  $$ \vert \trace (A^s -B^s) \vert \le \Vert M-N \Vert _1 
  \Vert \sum _{j=0}^{s-1}  A^{j} B^{s-1-j}\Vert.$$
  Or puisque $ T_N (\frac{1}{\tilde h_2})- T_N^{-1} \tilde h_2$ 
  est un oprateur positif nous avons 
  $$ \Vert A-B \Vert _1 \le \Vert T_N (h_1) \Vert 
  \trace \left(T_N^{-1} \tilde h_2 - T_N (\frac{1}{\tilde h_2}) \right)
  =o(N)$$
  d'aprs \cite{RS10}. D'autre part, toujours avec les m\^{e}mes 
  notations, 
  $$ \Vert A^j B^{s-1-j} \Vert \le\left( \Vert h_1\Vert_\infty 
  \Vert \frac{1}{\tilde h_2} \Vert _\infty\right)^s$$
  D'o finalement 
$$   \trace \left(T_N h_1 T_N ^{-1} \tilde h_2 \right)^s
-\trace \left( T_N h_1 T_N (\frac{1}{\tilde h_2} \right)^s
=o(N)$$
ce qui permet d'obtenir le rsultat avec le thorme \ref{thpuissance}
et en remplaant $\tilde h_2$ par sa valeur.

  \section{Dmonstration des thormes de grandes d\'eviations et de valeurs propres.}
  \subsection{d\'emonstration du lemme \ref{Szego}}
  Dans la suite nous poserons 
  $\Delta_{0} = ]-\delta_{0}^{-1} ,\delta_{0}^{-1} [$avec $\delta_{0} =2 K_{0} \Bigl \Vert {f_1}\Bigr \Vert _{\infty} \Bigl\Vert {f_2^{-1}} \Bigr\Vert_\infty$ la constante $K_{0}$ \'etant celle qui intervient dans l'\'enonc\'e du lemme \ref{norme1}.
  Dans la suite nous noterons $\lambda _{i}^{(N)}$, $1 \le i \le N$ les valeurs propres de $T_{N}(\frac{f_{1}}{f_{2}})$
   et $\mu_{i}^{(N)}$, $1 \le i \le N$ les valeurs propres de $T_{N}(f_{1}) T_{N}^{-1} (f_{2})$. 
   Nous pouvons alors \'enoncer le lemme 
    \begin{lemme} \label{fonda}
Pour tout $u\in \Delta_{0}$ $\alpha_{1}>0$ et $\alpha_{2}<0$  
  
$$ \sum_{l=1}^\infty \frac{(-1)^l }{l} u^l 
  \trace \left ( T_N f_1 T_N^{-1} f_2\right)^l = \sum_{l=1} ^\infty 
  \frac{(-1)^l }{l} u^l \trace T_N (\frac{f_1}{f_2}) ^l -
  \sum_{l=1} ^\infty \frac{1}{l} u^l \frac{ \Psi_1 ^{(l)}(0)}{l!}+o(1), $$
  \end{lemme}
  Dmontrons ce lemme .\\
 Soit $\epsilon>0$ t gr\^{a}ce  au lemme \ref{derivpuissance} et  puisque le rayon de convergence de la fonction $\Psi$ est sup\'erieur \`a $\delta_{0} $, nous savons qu'il existe $l_0$ tel que   $t\in \Delta_{0}$ et pour tout entier naturel $N$ nous avons les majorations 
 $$\Bigl \vert  \sum_{l=l_{0}}^\infty \frac{(-1)^l }{l} u^l 
 \left( \trace \left ( T_N f_1 T_N^{-1} f_2\right)^l - \trace T_N (\frac{f_1}{f_2}) ^l\right)\le \epsilon \quad 
 \mathrm{et} \quad  
  \sum_{l=l_{0}} ^\infty \frac{1}{l} u^l \frac{ \Psi_1 ^{(l)}(0)}{l!}\le \epsilon.$$
  D'autre part aprs le thorme \ref{thpuissance} il existe un $N_0$ tel que pour tout 
  $N \ge N_0$ et $s\le l_0$ on ait, 
  $$\Bigl \vert
 \Tr  \left( T_{N} (f_{1}) T_{N}^{-1} (f_{2}) \right)^s -
   \trace \left(T_{N} \left(\frac {f_{1}}{f_{2}}\right)^s\right)  - (-1)^{s+1}\frac {\Psi_{1}^{(s)} (0)}{s!} \Bigr \vert \le \frac{\epsilon}{l_0}.$$
  L'ingalit triangulaire permet alors de conclure.\\
  D'autre part pour tout $u \in \Delta_{0} $ nous pouvons \'ecrire 
  \begin{align*}
  \sum_{l=1}^\infty \frac{(-1)^l}{l} u^l \trace \left( T_{N}(\frac{f_{1}} {f_{2}})^l \right) &= \sum_{l=1}^\infty \frac{{-1}^l {l} }u^l \sum_{i=1} ^N
  (\lambda_{i}^{(N)})^l\\
  &= \sum_{i=1}^N  \sum_{l=0}^\infty \frac{{-1}^l }{l} u^l   (\lambda_{i}^{(N)})^l= \sum_{i=1}^N \ln (1_+ u \lambda_{i}^{(N)} ).
 \end{align*}
  De m\^{e}me nous obtenons 
  $$ \sum_{l=1}^\infty \frac{{-1}^l }{l} u^l \trace \left( T_{N} (f_{1}) T_{N}^{-1} (f_{2})\right)^l = \sum_{i=1}^N \ln (1+\mu_{i}^{(N)}).$$
  Nous avons donc obtenu 
  \begin{equation} \label{fonda2}
  \forall u \in \Delta_{0} \quad \lim_{N \rightarrow + \infty} \sum_{i=1}^N \frac{ \ln (1_+ u \lambda_{i)}^{(N)} )- \ln (1+ u \mu_{i}^{(N)})}{N}=0
  \end{equation}
  En reprenant \cite{GS} pp 62-63 nous pouvons conclure que pour tout fonction continue $F$ dans
  $ \Delta_{1}= [- \delta _{0}, \delta _{0}]$
  \begin{equation} \label {fonda3}
   \forall u \in \Delta_{1} \quad \lim_{N \rightarrow + \infty} \sum_{i=1}^N \frac{ F(u \lambda_{N)} )- F(\mu_{i}^{(N)})}{N}=0
\end{equation}
Ce qui traduit la convergence faible (ou en loi) de la suite de mesure $\sum_{i=1}^N \delta _{\lambda_{i}^{(N)}} $ vers la mesure 
image de la mesure de Lebesgue sur le tore par la fonction $ \frac{f_{1}}{f_{2}}$ (voir toujours \cite{GS} p65).
Nous pouvons alors appliquer la Proposition 3 pp79 de \cite{B.G.R.1} qui nous permet d'\'ecrire le th\'eor\`eme 
\ref{deviant1}
\subsection{D\'emonstration du th\'eor\`eme \ref{VP}}
Pour tout entier $N$ on suppose les valeurs propres $(\mu_{i}^{(N)})_{1\le i \le N}$ class\'ees par ordre croissante.
L'\'equation \ref{fonda3} peut encore 
s'\'ecrire 
 \begin{equation} \label {fonda4}
\forall u \in \Delta_{1} \quad \lim_{N \rightarrow + \infty} \sum_{i=1}^N \frac{ F(u \mu_{i}^{(N)})}{N} = \frac{1}{2 \pi} \int_{-\pi}^{\pi} F(f(x))dx
\end{equation}
  Ce qui implique 
  $$ \lim_{N\rightarrow \infty} \mu_{1}^{(N)} = 0, \quad \mathrm{et} \quad  \lim_{N\rightarrow \infty} \mu_{N}^{(N)} = \Vert 
  \frac{f_{1}}{f_{2}}\Vert_{\infty}.$$
  D\'emontrons cette propri\'et\'e. Posons $M = \Vert \frac{f_{1}}{f_{2}} \Vert_{\infty} $. On sait que les valeurs propres de $T_{N} (f_{1}) T_{N} ^{-1} (f_{2})$ sont dans $[0,M]$  pour cela on peut voir \cite{B.G.R.1} lemme 10 p 87 ou consid\'erer l'op\'erateur 
  $M T_{N } (f_{2}) -T_{N} (f_{1})$ qui est positif. Il en est de m\^{e}me pour l'op\'erateur 
  $M - T_{N}(f_{2})^{-1/2} T_{N} (f_{1}) T_{N}(f_{2})^{-1/2} $. Ses valeurs propres sont donc positives et on conclut en remarquant 
  que 
   $$ \mathrm{si} \quad T_{N} (f_{2})^{-1/2} T_{N} (f_{1}) T_{N}(f_{2})^{-1/2} (x) = \mu x \quad \mathrm{alors} \quad
   T_{N} (f_{1}) T_{N}(f_{2})^{-1} (y) =\mu y \quad \mathrm{avec} \quad
   x = T_{N}(f_{2})^{-1/2} (y).$$
   Supposons  maintenant que 
$\displaystyle{  \lim_{N\rightarrow \infty} \mu_{N}^{(N)}}\not=M $. 
Cela signifie qu'il existe un r\'eel $\epsilon>0$ et un entier $N_{0}$ tels qu'il existe 
  une sous-suite $\phi$ v\'erifiant ou 
  $$  \forall N \ge N_{0}Ê\quad  \mu_{\phi (N)}^{{\phi (N)}} - M < \epsilon.$$
  Une fonction $F$ continue, strictement positive, de support contenu dans $[M-\epsilon, M]$ nous donne alors une contradiction.
  On d\'emontre de m\^{e}me que $ \lim_{N\rightarrow \infty} \mu_{1}^{(N)} = 0$.
  \subsection {D\'emonstration du th\'eor\`eme \ref{deviant1}}
  Nous considrons donc la fonction 
  $ N L_N (t) = - \frac{1}{2} \sum_{i=1}^N 
  \ln (1 - 2 \mu_i^N t)$ o les $(\mu_i^N), \, 1\le N $ sont les valeurs propres de 
  $A_N =\left( T_N ^{1/2} f_1 T_N ^{-1}f_2 T_N ^{1/2} f_1\right)$ qui sont aussi 
  celles de $T_N f_1 T_N^{-1} f_2$. On crit, pour $\vert 2 t \vert < \frac{1}  {\max \vert \mu_{i}^N\vert} = \frac{1}{\delta }$
  \begin{align*} 
N L_{N}(t) = - \frac{1}{2} \sum_{i=1}^N \ln (1- 2 \mu_i ^N t ) &= \frac{1}{2}  \sum_{i=1} ^N \sum_{p=1} ^\infty 
  \frac{(2t)^p (\mu_i^N)^p}{p}  \\
  & =   \frac{1}{2} \sum_{p=1}^\infty \frac{(2t)^p }{p} 
  \trace \left ( T_N f_1 T_N^{-1} f_2\right)^p
  \end{align*} 
  En utilisant de nouveau le lemme \ref{fonda} il vient que pour tout $t \in \Delta_{0}$
     nous obtenons finalement 
  $$
  N L_{N}(t) =- \frac{1}{2} 
   \trace   T_{N} \ln \left( 1 - 2t \frac{f_1}{f_2} \right) -  \Psi (2t) +o(1).
$$
  Soit par analycit\'e pout tout $t\in \Delta$
  $$ L(t) = -\frac{1}{4 \pi} \int_{0}^{2 \pi} \ln \left( (1-t) - t \frac{f_{1}} {f_{2}} (\theta) \right) d \theta,$$
  et 
  $$
   \lim_{N\rightarrow + \infty } \left ( N L_{N}(t) + \frac{N}{4 \pi}\int_{0}^{2\pi} \ln \left( 1- 2 t \frac{f_{1}}{f_{2}} (\theta) \right) 
d\theta \right) = \frac{\Psi (2t)} {2}
  $$
  \section {Appendice}
 \subsection {Uniformit des restes dans la dmonstration du 
 lemme \ref{lemme1}}
  L'outil principal est l'valuation du reste de la formule d'Euler et
 Mac-Laurin que nous rappelons ici, pour deux entiers naturels
 $m$ et $n$ et une fonction $f$  
 \begin{align*}
 f(m) +f(m+1)+ \cdots +f(n) &= \int_m^n f(t) dt 
 + \frac{1}{2} \left(f(m)+f(n) \right) \\
 &+ \sum_{h=1}^r (-1) ^{h-1} \frac{B_h}{(2h)!} 
 \left( f^{(2h-1)} (n) - f^{(2h-1)} (m)\right) +R_r
 \end{align*}
 avec 
 $$ \vert R_r \vert \le \frac{2} {(2 \pi)^2}
  \int_m^n \vert f^{(r+1)} (t) \vert dt $$
  et o les $B_h$ sont les nombres de Bernouilli.\\
  \subsubsection{ Calcul de $A'_{1}$}
De part la nature des termes \`a consid\'erer le reste de la somme $A'_{1}$ est rapidement  traite gr\^{a}ce \`a la formule d'Euler et Mac-Laurin. 
Nous remarquons que pour avoir un reste en $o(N^{2\alpha_{2}-1}$ nous avons besoin de la condition 
 \begin{equation} \label{cond14}
 \delta > \frac{   \alpha_{2}}{\alpha_{2}-1}.
  \end{equation}

 \subsubsection {Reste du coefficient $A'_2$}

  Pour appliquer la formule d'Euler Mac-Laurin  au reste de la quantit 
  $A'_2$ on \'ecrit  
  \begin{align*} 
&  \sum_{l=0}^{N-s'} \beta _l ^{(\alpha_2)} 
  \beta _{l +s}^{(\alpha_2)} 
  \left( (1- \frac{l}{N} )^{\alpha_2}-1+ \alpha_2 \frac{l}{N} \right)
  =  \frac{N^{2\alpha_2-1}}{\Gamma^2(\alpha_2) c_2(1) }
\int_{0}^{m_0/N} t^{\alpha_2-1} 
 (t+s)^{\alpha_2-1} 
  \left( (1- t )^{\alpha_2}-1+ \alpha_2 t \right) dt \\
 & -  \sum_{l=0}^{m_0} \beta _l ^{(\alpha_2)} 
  \beta _{l +s}^{(\alpha_2)} 
  \left( (1- \frac{l}{N} )^{\alpha_2}-1+ \alpha_2 \frac{l}{N} \right) 
  + \frac{N^{2\alpha_2-1}}{\Gamma^2(\alpha_2) c_2(1) }
 \int _0 ^{1-\frac{s'}{N}} t^{\alpha_2 -1} 
 (t+ \frac{s}{N} )^{\alpha_2-1} 
 \left( (1-t)^{\alpha_2}-1+\alpha_2 t\right)dt \\
 &+ m_0 ^{\alpha_2 -1} 
 (m_0 + s )^{\alpha_2-1} 
 \left( (1-\frac{m_0}{N})^{\alpha_2}-1+\alpha_2 \frac{m_0}{N}\right)\\
 &+ (N-s') ^{\alpha_2 -1} 
 ((N-s')+ s )^{\alpha_2-1} 
 \left( (1-\frac{N-s'}{N})^{\alpha_2}-1+\alpha_2 \frac{N-s'}{N}\right)
  +R_0 
  \end{align*}
  o\`u $R_{0}$ est le reste d'Euler et Mac-Laurin, d\'efini comme ci-dessus.
  Pour rendre les encadrements des diff\'erentes quantit\'es plus ais\'es on divise l'ensemble des indices $s$ 
 en deux intervalles $[0, N\delta_1]$ et 
 $[N \delta_1, N - N ^\delta]$ avec $0<\delta <1$ et $0<\delta _{1}<1$. La quantit\'e  $N^\delta $ correspond \`a l'entier $n_{0}$ de la d\'emonstration.  Nous ferons 
\'egalement intervenir le r\'eel $m_{0} =N^{\delta _{0}}$ ($0<\delta _{1}<1$) qui soit tel 
que pour $m\ge m_{0}$ l'usage de l'asymptotique de $ \beta_{m}^{(\alpha)}$ soit pertinent. 
Nous avons bien s\^{u}r 
$0<\delta _{0}<1$, $0<\delta _{1}<1$, $N^{\delta _{0}}<N \delta_{1} <N-N^\delta $.
  La quantit\'e $\delta $ est commune aux quatre calculs des restes relatifs aux termes 
  $ \quad A'_{2}, \quad A'_{3}, \quad A'_{4}$. Par contre les quantit\'es $\delta_{0}$ et $\delta_{1}$
  peuvent \^{e}tre choisies ind\'ependamment dans chacun des quatre calculs.    
  Les majorations des diff\'erents termes intervenants dans la somme sont classiques et on obtient l'approximation annonc\'e\'e avec un reste d'ordre 
  $o(N^{2\alpha_{2}-1})$
  uniform\'ement par rapport \`a $s$.
    Pour terminer l'\'etude du reste coefficient $A'_{2}$ il nous faut encore consid\'erer pour $s \in [0, N-N^\delta]$ le reste de 
 l'\'egalit\'e 
 \begin{align*}
 \frac{1}{N} \sum _{N-s'+1}^{+ \infty} \beta_{l}^{(\alpha_{2})} \overline{\beta_{l+s}^{(\alpha_{2})}} l &=
 \frac{1}{N} \int _{N+s'+1}^{+ \infty} \frac{t^{\alpha_{2}} (t+s)^{\alpha_{2}-1}} {\Gamma^2(\alpha_{2}) c_{2}(1)} dt \\
 & + \frac{1}{2N} \frac{ (N-s')^{\alpha_{2}} (N+s-s')^{\alpha_{2}-1} }{\Gamma^2 (\alpha_{2}) c_{2}(1)}+R_{2}.
 \end{align*}
 Nous avons imm\'ediatement 
 \begin{align*}
  \frac{1}{2N} \frac{ (N-s')^{\alpha_{2}} (N+s-s')^{\alpha_{2}-1} }{\Gamma^2 (\alpha_{2}) c_{2}(1)} 
 & = O \left( N^{\delta _{0} \alpha_{2}} N^{\alpha_{2}-2} \right) =o(N^{2\alpha_{2}-1})\; \mathrm {si} \; \alpha_{2}<0\\
 &= O(N^{2\alpha_{2} -2}) =o(N^{2\alpha_{2}-1})\; \mathrm {si} \; \alpha_{2}>0.
\end{align*}
 et d'autre part 
\begin{equation} \label{restedesrestes1}
 \vert R_{2} \vert \le \frac{1}{N} \int_{N-s'+1}^{+ \infty} \frac{1}{\Gamma^2 (\alpha_{2})c_{2}(1) } 
 \left( \alpha_{2} t^{\alpha_{2}-1} (t+ \frac{s}{N})^{\alpha_{2}-1} + (\alpha_{2}-1) t^{\alpha_{2}} (t+s)^{\alpha_{2}-2} \right) dt. 
\end{equation}
La d\'eriv\'ee intervenant dans l'\'equation \ref{restedesrestes1} \'etant de signe constant nous avons 
 \begin{align*}
 \vert R_{2} \vert &\le  \frac{1}{2N} \frac{ (N-s')^{\alpha_{2}} (N+s-s')^{\alpha_{2}-1} }{\Gamma^2_{\alpha_{2}} c_{2}(1)} \\
 &= O \left( N^{\delta _{0} \alpha_{2}} N^{\alpha_{2}-2} \right) =o(N^{2\alpha_{2}-1})\; \mathrm {si} \; \alpha_{2}<0\\
 &= O(N^{2\alpha_{2}-2}) =o(N^{2\alpha_{2}-1})\; \mathrm {si} \; \alpha_{2}>0.
\end{align*}

  \subsubsection {Reste du coefficient $A'_3$}
Ici aussi pour faciliter nos calculs nous sommes oblig\'e de diviser  l'intervalle auquel appartient le
param\`etre $s$ en $[0, N \delta _{1}]$, et en $ [N \delta _{1}, N-N^\delta]$.
 Nous cherchons tout d'abord \`a majorer uniform\'ement le reste $Q$ dans la formule 
\begin{align*}
&\sum_{l=0}^{N-s'} \beta_{l}^{(\alpha_{2})} \beta_{l+s}^{(\alpha_{2})} 
\left( (1- \frac{l+s}{N})^{\alpha_{2}} - (1-\frac{s}{N})^{\alpha_{2}} +\alpha_{2} \frac{l}{N}  (1-\frac{s}{N} )^{\alpha_{2}-1}\right)=\\
&= \frac{N^{2 \alpha_{2}-1}}{c_{2}(1) \Gamma^2(\alpha_{2})} 
\int_{0}^{1-\frac{s'}{N}} t^{\alpha_{2}-1} ( t+ \frac{s}{N})^{\alpha_{2}-1} 
\left( (1-t-\frac{s}{N})^{\alpha_{2}} -(1-\frac{s}{N})^{\alpha_{2}}\right.+\\
 &+\left. \alpha_{2} t (1-\frac{s}{N})^{\alpha_{2}-1}  \right) dt
 +Q.
 \end{align*}
 Pour ce faire nous utilisons le m\^{e}me type de dcomposition 
 que dans le cas du reste de $A'_{2}$ avec $m_0=N^{\delta _{0}}$  
 un entier fix comme ci-dessus.
 Comme pr\'ec\'edemment nous sommes conduits \`a utiliser la formule d'Euler et Mac-Laurin et nous \'ecrivons 
\begin{align*}
Q &= \frac{m_{0}^{\alpha_{2}-1}}{c_{2}(1) \Gamma^2(\alpha_{2})}  
 (m_{0}+s)^{\alpha_{2}-1} 
 \left( (1-\frac{m_{0}+s}{N} )^{\alpha_{2}} - (1-\frac{s}{N})^{\alpha_{2}} +
 \frac{\alpha m_{0}} {N} (1-\frac{s}{N})^{\alpha_{2}-1} \right) \\
& +   \frac{(N-s')^{\alpha_{2}-1}}{c_{2}(1) \Gamma^2(\alpha_{2})}  
 (N-s'+s)^{\alpha_{2}-1} 
 \left( (1-\frac{N-s'+s}{N} )^{\alpha_{2}} - (1-\frac{s}{N})^{\alpha_{2}} +
 \frac{\alpha N-s'} {N} (1-\frac{s}{N})^{\alpha_{2}-1} \right) \\
& +  \sum_{l=0}^{n_{0}} \beta_{l}^{(\alpha_{2})} \beta_{l+s}^{(\alpha_{2})} 
  \left( (1- \frac{l+s}{N})^{\alpha_{2}} -(1-\frac{s}{N})^{\alpha_{2}} +\alpha_{2} \frac{l}{N} ( 1- \frac{s}{N} ){\alpha_{2}-1}
  \right) ds\\
  &- \frac{N^{2\alpha_{2}-1} } {c_{2}(1) \Gamma^2 (\alpha_{2}) } 
  \int_{0} ^{n_{0}/N} t^{\alpha_{2}-1} (t+\frac{s}{N} )^{\alpha_{2}-1} 
  \left( (1-t-\frac{s}{N} )^{\alpha_2} -( 1- \frac{s}{N})^{\alpha_{2}}+\alpha_{2} t (1-\frac {s}{N})^{\alpha_{2}-1}\right)
  dt\\
  +R'_{0}
  \end{align*}
  avec $R'_{0}$ qui correspond au reste de la formule d'Euler et Mac-Laurin. 
   \begin{align*}
R'_{0}&\le \frac{N^{2\alpha_{2}-2} } {c_{2}(1) \Gamma^2 (\alpha_{2}) }  \int_{m_{0}/N} ^{1-s'/N} 
 \Big \vert \left( (\alpha_{2}-1) t^{\alpha_{2}-2} (t+ \frac{s}{N})^{\alpha_{2} -1} (2 t+\frac{s}{N}) \right) \\
& \left ( (1-t-\frac{s}{N})^{\alpha_{2}} -(1-\frac{s}{N})^{\alpha_{2}}
 + \alpha_{2} t (1-\frac{s}{N})^{\alpha_{2}-1}  \right) dt \Big \vert  \\
&+  \int_{m_{0}/N} ^{1-s'/N} \Big \vert
\left( \alpha_{2} t^{\alpha_{2}-1} (t+ \frac{s}{N})^{\alpha_{2} -1} \right) \\
& \left( -\alpha_{2} (1-t-\frac{s}{N})^{\alpha_{2}-1} 
 + \alpha_{2} (1-\frac{s}{N})^{\alpha_{2}-1}  \right) \Big \vert dt 
 \end{align*}

  La majoration de la quantit\'e $\vert R'_{0} -Q$ et classique et nous donne les conditions 
  \begin{equation} \label{cond6}
  \delta > \frac{2\alpha_{2}} {\alpha_{2}-2}
\quad 
  \delta > \frac{2\alpha_{2}} {(\alpha_{2}-1)}
  \end{equation}
qui sont compatibles entre elles et avec la condition \ref{cond14}. Toutes ces conditions nous donne une valeur 
plus pr\'ecise de $N^\delta $ la d\'efinition de $N^\delta $.
 Pour majorer  uniform\'ement le reste de la formule d'Euler et Mac-Laurin pour $s\in [N\delta_{1} , N^\delta ]$ on \'ecrit 
 $ R'_{0} \le \frac{N^{2\alpha_{2} -2}}{c_{2}(1) \Gamma ^2 (\alpha_{2})} (\vert I_{1}+ I_{2}).$ 
 avec 
 $$ I_{1} =   \int_{m_{0}/N} ^{1-s'/N} 
 \Big \vert \left( (\alpha_{2}-1) t^{\alpha_{2}-2} (t+ \frac{s}{N})^{\alpha_{2} -1} (2 t+\frac{s}{N}) \right) 
 \left ( (1-t-\frac{s}{N})^{\alpha_{2}} -(1-\frac{s}{N})^{\alpha_{2}}
 + \alpha_{2} t (1-\frac{s}{N})^{\alpha_{2}-1}  \right)  \Big \vert dt$$
 et 
 $$ I_{2 }= 
 \int_{m_{0}/N} ^{1-s'/N} \Big \vert
\left( \alpha_{2} t^{\alpha_{2}-1} (t+ \frac{s}{N})^{\alpha_{2} -1} \right) 
 \left( -\alpha_{2} (1-t-\frac{s}{N})^{\alpha_{2}-1} 
 + \alpha_{2} (1-\frac{s}{N})^{\alpha_{2}-1}  \right) \Big \vert dt 
 $$
 Si $\delta_u$ dsigne le coefficient d'indice $u$ du dveloppement en srie entire de $(1-v)^{\alpha_2}$ 
 nous pouvons crire 
$$
 (1-t-\frac{s}{N})^{\alpha_{2}} -(1-\frac{s}{N})^{\alpha_{2}}
 + \alpha_{2} t (1-\frac{s}{N})^{\alpha_{2}-1} 
 = \sum_{n \ge 0} \delta_{n+2} t^{n+2} (1-\frac{s}{N})^{\alpha_2-n}
 $$
 et d'autre part si $t \in [\frac{m_0}{N}, 1-\frac{s'}{N}]$ nous avons
 la majoration 
 $$
 I_{1} \le O \left ( \int_{m_{0}/N} ^{1-s'/N}  \sum_{n\ge 0}\delta _{n} t^{\alpha_{2}+n} (1-\frac{s}{N})^{\alpha_{2}-n} dt \right)$$
 ou encore 
 $$
 I_{1}=
  O \left( (1-\frac{s}{N})^{2\alpha_{2}+1}\right) + O \left( \sum _{n \ge 0}   (\frac{m_{0}} {N}) ^{\alpha_{2}+n+1} (\frac{N^\delta }{N})^{\alpha_{2}-n}\right) = o(1)
 $$
 si l'on choisit $\delta $ tel que $N^\delta >m_{0}$.
  Nous avons de m\^{e}me 
 $$
 -\alpha_{2} (1-t-\frac{s}{N})^{\alpha_{2}-1} 
 + \alpha_{2} (1-\frac{s}{N})^{\alpha_{2}-1}
 = - \left( \sum_{n \ge 0} \delta _{n+1} t^{n+1} 
 (1-\frac{s}{N})^{\alpha_2 -n-1} \right)
 $$
 ce qui donne 
 $$
 I_{2}= O \left( \int_{m_{0}/N} ^{1-s'/N}  \sum_{n\ge 0}\delta _{n} t^{\alpha_{2}+n} (1-\frac{s}{N})^{\alpha_{2}-n-1} dt \right) =o(1)$$
 
 comme pour $I_{1}$.
\subsubsection{ Calcul de $A'_{4}$}

 Nous devons utiliser la dcomposition 
\begin{align*}
&\sum_{l=0}^{N-s'} \beta_l^{(\alpha_2)} 
\overline{ \beta_{l+s} ^{(\alpha_2)}} 
\left( (1-\frac{l+s}{N} )^{\alpha_2} -
(1-\frac{s}{N})^{\alpha_2} \right) =\\
&= \sum_{l=0}^{m_0} \beta_l^{(\alpha_2)} 
\overline{ \beta_{l+s} ^{(\alpha_2)}} 
\left( (1-\frac{l+s}{N} )^{\alpha_2} -
(1-\frac{s}{N})^{\alpha_2} \right) \left( (1-\frac{l}{N})^{\alpha_2}-1\right)\\
& - \int_0^{m_0} l^{\alpha_2-1} (l+s)^{\alpha_2-1}
\left( (1-(\frac{l+s}{N} )^{\alpha_2} - 
(1-\frac{s}{N})^{\alpha_2} \right)\left( (1-\frac{l}{N})^{\alpha_2}-1\right)dl
 \\
&+  m_0^{\alpha_2} 
(m_0+s) ^{\alpha_2} 
\left( (1-\frac{m_0+s}{N} )^{\alpha_2} -
(1-\frac{s}{N})^{\alpha_2}\right)\left( (1-\frac{m_0}{N})^{\alpha_2}-1\right)\\
&+ (N-s')^{(\alpha_2} 
(N-s'+s) ^{\alpha_2} 
\left( (1-\frac{N-s'+s}{N} )^{\alpha_2} -
(1-\frac{s}{N})^{\alpha_2}\right)
\left( (1-\frac{N-s'}{N})^{\alpha_2}-1\right)\\
&+R
\end{align*}
o 
$R$ est le reste d'Euler et Mac-Laurin.
\\
Comme d'habitude nous distinguons le cas 
$s \in [0, N\delta_1] $ et $s\in [N \delta_1, N-N^\delta]$. 
Nous obtenons la condition 
 \begin{equation} \label{cond16}
 \delta > \frac{2\alpha_2}{2 \alpha_2 -1}.
\end{equation}
     \bibliography{Toeplitzdeux}

 \end{document}